\documentclass[11pt,fleqn,a4paper]{article} 
\usepackage{amsfonts,amsmath}
\usepackage[english]{babel}
\usepackage[T1]{fontenc}
\usepackage{dsfont}
\usepackage[latin1]{inputenc}
\usepackage{a4}
\usepackage{mathbbol}
\usepackage{amsfonts,amssymb,amsmath,epsfig}
\usepackage{bbold}
\usepackage{stmaryrd,epsfig}
\usepackage{pifont,fancyheadings,graphicx,subfigure,ifthen,epsfig}
\usepackage{enumerate,rotating,hangcaption} 
\newcommand{\f}{\mathbf}

\newcommand{\ep}{\epsilon}

\newcommand{\Dx}{{\partial_x}} 
\newcommand{\Dy}{{\partial_y}} 
\newcommand{\B}{\bar} 
\newcommand{\Dxx}{{\partial^2_{xx}}} 
\newcommand{\Dyy}{{\partial^2_{yy}}} 
\newcommand{\Dxy}{{\partial^2_{xy}}} 
\newcommand{\dx}{{\partial_x}} 
\newcommand{\dy}{{\partial_y}} 
 
\newcommand{\3}{{ {\mathds 1}}_{x \geq 0}} 
\newcommand{\5}{{ {\mathds 1}}_{x \leq 0}} 
\newcommand{\4}{{ {\mathds 1}}_{y \geq 0}}

\newcommand{\R}{\mathds R} 
\newcommand{\C}{\mathds C}

\newtheorem{proposition}{Proposition} 
 
\newtheorem{coroll}{Corollary}
\newenvironment{proof}{{\bf Proof.}}{$\diamondsuit$} 

\begin{document}

 
\title{Simulation of Laser Beam Propagation With a Paraxial Model in a Tilted Frame}

\author{Marie Doumic \thanks{Corresponding author M. Doumic: marie.doumic@inria.fr} \thanks{I.N.R.I.A Rocquencourt, BANG project, Domaine de Voluceau, B.P. 105, 78153 Rocquencourt, France, marie.doumic@inria.fr} \qquad Fr\'ed\'eric Duboc \thanks{CEA/Bruy\`eres, B.P. 12, 91680 Bruy\`eres-Le-Chatel, France, frederic.duboc@cea.fr} \qquad
Fran\c{c}ois Golse \thanks{Ecole Polytechnique, Centre de Math\'ematiques Laurent Schwartz, 91128 Palaiseau 
Cedex, France, golse@math.polytechnique.fr} \qquad
R\'emi Sentis \thanks{CEA/Bruy\`eres, B.P. 12, 91680 Bruy\`eres-Le-Chatel, France, remi.sentis@cea.fr}}

\maketitle 

\begin{abstract}
We study  the Schr\"odinger equation which comes from  the paraxial approximation of 
the Helmholtz equation in the case where the direction of propagation is tilted with 
respect to the boundary of the domain.  In a first part, a mathematical analysis is made which leads to an analytical formula of the solution  in the simple case where the refraction index and the absorption coefficients are constant. Afterwards, we propose a numerical 
method for solving the initial problem which uses the previous analytical expression. Numerical results are presented. We also sketch an extension to a time dependant model which is relevant for laser plasma interaction.
\end{abstract}







                                \section{Introduction}


For the simulation of the propagation of a monochromatic laser beam in a
medium where the local refractive index is nearby a constant, it is
classical to use the paraxial approximation of the Maxwell equations. This
approximation takes into account diffraction and refraction phenomena ; it
is intensively used for decades in optics and in a lot of models related to
laser-plasma interaction in Inertial Confinement Fusion experiments (cf \cite{para-3},\cite{dorr}, \cite{ria}, \cite{feu} and the bibliography of these
references). Let us first recall briefly the outlines of this approximation.
Denote by $2 \pi \epsilon $ 
the laser wave-length, it is in the
order of 1 $\mu m$ and is very small compared to the characteristic length
of the simulation domain (which is in the order of some $mm$ for the
Inertial Confinement plasmas). According to laws of optics, the laser
electromagnetic field may be modeled by the solution $\psi $ of the
following Helmholtz equation (which comes from the time envelope of the full
Maxwell equations): 
\begin{equation}
\epsilon ^{2}{\Delta }\psi +\psi +2i\epsilon \nu _{t}\psi =0,  \label{base}
\end{equation}%
where we have denoted:
\begin{equation*}
\nu _{t}(\mathbf{x})=\nu (\mathbf{x})+i\mu (\mathbf{x}),
\end{equation*}%
so $\nu_t$ is a complex function, its real part $\nu $ corresponds to a conveniently
scaled absorption coefficient and its imaginary part $\mu $ to the variation
of the refractive index ($1-2\epsilon \mu $ is equal to the square of the
refractive index $n$ up to a multiplicative constant).

We assume also that the light propagates according a fixed direction defined
by the unit vector $\mathbf{k}.$ After making the classical WKB expansion: 
\begin{equation}
\psi =u\exp (\frac{i\mathbf{k}.\mathbf{x}}{\epsilon }),  \label{aaa}
\end{equation}%
equation (\ref{base}) may read as $2i\nu _{t}u+2i\mathbf{k}.{\mathbf{\nabla }%
}u+\epsilon \Delta _{\bot }u=\epsilon (\mathbf{k}.{\mathbf{\nabla }})^{2}u,$
where $\Delta _{\bot }$ is the Laplace operator with respect to the
transverse variable: 
\begin{equation*}
\Delta _{\bot }\bullet =\nabla .[(\mathbf{1}-\mathbf{k}\otimes \mathbf{k}%
)\nabla \bullet ],\quad \mathbf{1}\hbox{ being the unit
 diagonal  tensor.}
\end{equation*}%
Assuming that $u$ is slowly varying with respect to the longitudinal
variable, we can neglect the right hand side of the previous equation.
Therefore $u$ satisfies the classical paraxial equation for wave propagation:
\begin{equation}
i\mathbf{k}.\nabla u+\frac{\epsilon }{2}\Delta _{\bot }u+i\nu _{t}u=0, \qquad {\rm with }  \quad  \nu _{t}=\nu +i \mu.
\label{bbb}
\end{equation}

For this kind of model, it is usual to handle a simulation box which is a
parallelepiped and the laser beam is assumed to enter into the simulation box on
a plane boundary denoted by $\Gamma _{0}.$ Let us denote $\mathbf{n}$ the
outward normal vector to the incoming boundary $\Gamma _{0}.$ Classically,
the crucial assumption is that the laser beam enters into  the simulation domain
with a very small incidence angle, that is to say the vector $\mathbf{k}$ is
almost equal to $-\mathbf{n}$. Then, in such a framework, (\ref{bbb}) is a
classical linear Schr\"{o}dinger equation, the operator $\mathbf{k}.\nabla $
plays the part of time derivative and the boundary condition on $\Gamma _{0}$
which reads $u=u^{in}\ $ (where $u^{in}$ is a given function defined on $%
\Gamma _{0}$) plays the part of the initial condition. On the other hand,
artificial absorbing boundary conditions are to be imposed on the faces of
the simulation domain parallel to the vector $\mathbf{k},$ (see for example 
\cite{Arnold-2}, \cite{DM}, \cite{had}). The numerical methods are always
implemented on an orthogonal mesh and are based on a splitting with respect
to the main spatial variable between the diffraction part $(\frac{\epsilon }{%
2}\Delta _{\bot }u)$ and refraction part $(i\nu _{t}u),$ see \cite{para-3}, 
\cite{fcs}, \cite{dorr} for example.

We address in this paper a different case where the incidence angle of $%
\mathbf{k}$ with $-\mathbf{n}$ is large; these simulations are called tilted
frame simulations. This kind of simulations is of particular interest if one
has to deal with the crossing between two beams (in the high energy laser
devices, a large number of beams are focused on the target, therefore beam
crossing may be taken into account, see \cite{lin} for a survey on
related laser propagation problems); an example of such simulations in a
very simplified case may be found on Figure \ref{beau}. This tilted frame model has
been considered some years ago by physicists for dealing with beam crossing
problems (see \cite{tik}).

Simulations in a tilted frame are also necessary for dealing with special
situations. For instance for the propagation of a beam in a domain where
the profile of the refractive index $n$ is such that $n^{2}(\mathbf{x}%
)=n_{0}^{2}(1-\varepsilon \mu (\mathbf{x}))$ (with $n_{0}$ constant smaller
than 1) in a first subdomain $\mathcal{D}$ and $n^{2}(\mathbf{x})=\mathcal{N}%
(\mathbf{x}.\mathbf{n}^{\ast })+\delta \mathcal{N}(\mathbf{x})$ (where $%
\mathcal{N}\in \lbrack 0,n_{0}]$ depends on a one-dimension variable $%
\mathbf{x}.\mathbf{n}^{\ast }$ and $\delta \mathcal{N}$ is small with
respect to $1$) in a second juxtaposed subdomain $\mathcal{D}^{H}$, one must
handle the paraxial equation (\ref{bbb}) in subdomain $\mathcal{D}$ and the
Helmholtz equation (\ref{base}) in subdomain $\mathcal{D}^{H}$. For the
numerical solution of (\ref{base}), one has to solve a huge linear system
(corresponding to the discretization of the equation on a very fine grid)
and for handling this huge linear system, it is necessary that the variable $%
\mathbf{x}.\mathbf{n}^{\ast }$ corresponds to one of the main direction of
$\mathcal{D}^{H}$. Therefore the full simulation on $(\mathcal{D}%
\cup \mathcal{D}^{H})$ has to be performed in a box such that the corresponding normal vector $\mathbf{n}$ must be parallel to $\mathbf{n}^{\ast }$   (see \cite{desr}
for details for this kind of simulations).

In the case of a large incidence angle, the crude expansion $\psi =U\exp ({-i%
\mathbf{n}.\mathbf{x}}/{\epsilon })$ leads to difficulties and to overcome
these difficulties, it has been proposed in \cite{faf} to replace the
transverse Laplacian by a pseudodifferential operator, but with this
approximation, $U$ is not slowly varying with respect to the spatial
coordinates therefore it is necessary to handle very fine mesh -at least 10
cells per wave length- to get accurate results. One can also refer to the
works in the spirit of \cite{lee} in the acoustic framework but the
application to the optics problems seems to be difficult.

Here we consider the expansion $\psi =u\exp ({i\mathbf{k}.\mathbf{x}}/{%
\epsilon }),$ with $u$ slowly varying with respect to $\mathbf{k}.\mathbf{x,}
$ so we have to deal with the tilted frame Laplace operator $\Delta _{\bot }$
and one has to supplement the equation (\ref{bbb}) with a right incoming
boundary condition on $\Gamma _{0}$. For the statement of this boundary
condition, one assumes that a fixed plane wave $\psi ^{in}=u^{in}\exp (i%
\mathbf{k}.\mathbf{x}/{\epsilon })$ enters into the domain where $u^{in}$
is a given function of the variable which is orthogonal to $\mathbf{k}.$ Now, for the
Helmholtz problem, the boundary condition is classical and may be written as 
$(\epsilon \mathbf{n}.\nabla +i\mathbf{k}.\mathbf{n})(\psi -u^{in}e^{i%
\mathbf{k}.\mathbf{x}\mathbf{/\epsilon }})=0,$ then  using (\ref{aaa}) and
an asymptotic expansion with respect to the small parameter ${\epsilon }$,
the corresponding boundary condition for equation (\ref{bbb}) may read in a
natural way as: 
\begin{equation}
(\epsilon \mathbf{n}.\nabla _{\bot }+2i\mathbf{k}.\mathbf{n})(u-u^{in})=0,
\label{Cbord}
\end{equation}

where $\nabla _{\bot }=\nabla -\mathbf{k}(\mathbf{k}.\nabla )$ denotes the
gradient orthogonal to $\mathbf{k}.$ See \cite{D} for a justification of the
paraxial approximation in the special case we are dealing with.

If one sets $\mathbf{x}=(x,y,z)$ in 3D and $\mathbf{x}=(x,y)$ in 2D, the
entrance boundary $\Gamma _{0}$ corresponds in this paper to $x=0$. In the
sequel we consider a 2D problem but most of the ideas of this work may be
extended to the 3D case.

Equation (\ref{bbb}) may be recast as:
\begin{equation*}
i(k_{x}\partial _{x}u+k_{y}\partial _{y}u)+\frac{\epsilon }{2}\Delta _{\bot
}u+i\nu _{t}u=0,
\end{equation*}%
and up to our knowledge, the numerical solution of this kind of equations is
novel; the main difficulty is to handle correctly the tilted Laplace
operator $\Delta _{\bot }u$. For the mathematical analysis of the problem,
one key result is the following (cf. proposition 2). On the half-space $%
\{(x,y)\quad \mathrm{s.t.}\quad x\geq 0\},$ if the coefficient $\nu _{t}$ is
a positive real constant, after taking the Fourier transform with respect to
the $y$ variable, the problem~(\ref{bbb})(\ref{Cbord}) is equivalent to an
ordinary differential equation with respect to the $x$ variable and it is
possible to exhibit an analytical solution. This analytical formula is the
convenient tool for numerical treatment of the diffraction part of (\ref{bbb}%
) in the general case where $\nu _{t}$ is not constant.

The paper is organized as follows. In Section \ref{section-demi-2}, after
setting classical energy estimates for Problem~(\ref{bbb}) supplemented
by (\ref{Cbord}), we prove the above mentioned theoretical result. 

 Section~\ref{section-num} is devoted to
the description of the numerical scheme for solving Problem (\ref{bbb})(%
\ref{Cbord}) ; it is based on a splitting method with respect to the spatial
variable $x$ using fast Fourier transforms on a first step (for the
diffraction part) and a standard finite difference method on a second step
(for the advection and refraction part).

In Section \ref{sec:tilted}, we give the numerical results on the initial problem and for a
model where the coefficient $\mu $ in (\ref{bbb}) is replaced by $f(|u|)$
corresponding to the autofocusing which occurs in the laser-plasma
interaction (see \cite{rose} for instance). From a physical point of view,
this term represents a variation of the plasma electronic density caused by
the ponderomotrice force of the laser. In the last section we consider a
more general model where the stationary problem (\ref{bbb}) is replaced by a
time dependent one which is coupled to a hydrodynamic system for a suitable
modeling of the plasma behavior.

                                 \section{Analysis of the Tilted
				 Paraxial Equation}

\label{section-demi-2}
For reasons which will appear in the sequel, we  assume in
this section that
\begin{equation}\label{hyp1}
\hbox {inf}_{\f{x}} \nu(\f{x})  > 0.
\end{equation}
 We first study the problem where the simulation domain is the half-space: 
$$
{\cal D}=\big\{  \f{x}  =(x,y) \quad  {\rm s.t.} \quad  x>0 \big\},\quad \Gamma_0=\big\{\f{x}=(0,y)\big\}.
$$
Assuming that $\mu$ is a bounded function,  we consider the  following problem:
\begin{eqnarray}
i\f{k}\cdot\f{\mathbf{\nabla }}u+\frac{\epsilon }{2}%
\Delta _{\bot }u- \mu u+i\nu u =0  \qquad   {\rm on} \;  {\cal D},
\label{bbbb}
\\
\label{CE0}
(i\epsilon \f{n}. \nabla _{\bot } -2 \f{k}.\f{n} )
(u-u^{in})=0  \qquad   {\rm on} \; \Gamma_0.
\end{eqnarray}

                \subsection{Energy Estimate}
Let us first state the following classical estimate.

\label{energie-demi}

\begin{proposition}\label{energie-demi-theo}
Let $(i\ep \f{n}. \nabla _{\bot } - 2\f{k}.\f{n})u^{in} \in L^2 (\R)$. 
If  $u \in H^1({\cal D})$ is a solution  to Problem (\ref{bbbb}) (\ref{CE0}), it is unique. Moreover, we have the following
 stability estimate, with a constant $C$ independent of $\nu,$  $\mu$:
$$ \iint\limits_{\mathcal{D}}2\nu |u |^{2} + \int\limits_{\Gamma_0}  |\f{k}\cdot\f{n}| |u|^2 dy \leq C \int \limits_{\Gamma_0} |(i\ep \f{n}. \nabla _{\bot }
 - 2\f{k}.\f{n}) u^{in}|^2 dy.$$
\end{proposition}

\begin{proof}
Let us denote $D=\f{n}. \nabla _{\bot }$.  
Doing the scalar product of Equation (\ref{bbb}) with $u$ and taking  its imaginary part, we get:
$$ \int\limits_{\Gamma_0}  \biggl(|u|^2 \f{k} \cdot \f{n}+ \frac{\ep}{ 2i} (\B{u} D u - u D \B{u})\biggr)dy +  \iint\limits_{\cal D} 2 \nu |u|^2 d\f{x} =0.$$
According to the boundary condition (\ref{CE0}) we check that:
$$ \frac{\ep}{2i} (\B{u} D u - u D \B{u}) =-2 \f{k}\cdot \f{n} |u|^2 + {\cal I}m \bigl(\B{u}(\ep D + 2i \f{k}\cdot \f{n})u^{in}\bigr).$$
Then we get:
\begin{equation}
 \iint\limits_{\mathcal{D}}2\nu |u |^{2} d\f{x} +
 \int\limits_{\Gamma_0}  |\f{k}.\f{n}| |u|^2  dy = -{\cal I}m \bigl( \int\limits_{\Gamma_0} \B{u}
 (\ep D + 2i\f{k}.\f{n}){u^{in}} dy \bigr).
\label {sec}
\end{equation}
According to (\ref{sec}), if $(i\ep D -
2\f{k}\cdot\f{n})u^{in} =0,$ we see that $\iint\limits_{\mathcal{D}}2\nu |u
|^{2}d\f{x}=0,$ so $u=0.$ Therefore we get the uniqueness of the solution of Problem (\ref{bbbb})(\ref{CE0}).
 
To obtain the 
stability inequality, we first see that Equation (\ref{sec}) implies:
$$|\f{k}\cdot\f{n}|  \int\limits_{\Gamma_0}   |u|^2  \leq \sqrt{ \int\limits_{\Gamma_0} |u|^2 } \sqrt{\int\limits_{\Gamma_0}|(\ep D +
 2i\f{k}\cdot\f{n}){u^{in}}|^2 }.$$ 
Using this estimate, Equation (\ref{sec}) leads to:
$$\iint\limits_{\mathcal{D}}2\nu |u |^{2} d\f{x} +
\int \limits_{\Gamma_0}  |\f{k}\cdot\f{n}| |u|^2   \leq  \sqrt{ \int\limits_{\Gamma_0} 
 |u|^2 }\sqrt{\int\limits_{\Gamma_0}  |(\ep D +
 2i\f{k}\cdot\f{n}){u^{in}}|^2 } 
\leq \frac{1}{ |\f{k}\cdot \f{n}|} \int\limits_{\Gamma_0}  |(\ep D + 2i\f{k}\cdot\f{n}){u^{in}}|^2 .
$$ 
\nopagebreak
\end{proof}

By the same technique we get also the following estimate:
\begin{equation*}
\iint\limits_{\mathcal{D}}2\nu |u |^{2}
+\int\limits_{\Gamma_0}\frac{|\f{k}.\f{n}|}{2}|  
\frac{(i\epsilon D + 2 \f{k}.\f{n} )u}{2|\f{k}.\f{n}|}| ^2
 = 
\int\limits_{\Gamma_0 }  |\f{k}.\f{n}|\biggl(|u|^2 
+\frac{1}{2}  
|\frac{(i\epsilon D  - 2 \f{k}.\f{n}) u^{in} }{2|\f{k}.\f{n}|}| ^2
\biggr),
\end{equation*}
which says that the absorbing energy plus the the outgoing energy is equal to the incoming energy.

\subsection{ Analytical Form of the Solution in the Case $\nu_t$ Constant}

\label{subsection-fourier-demi}
 We now assume that $\mu=0$ and  $\nu $ is constant for getting an analytical form of the solution to Problem~(\ref{bbb})(\ref{Cbord}).
We denote  $\f{k}=(k_x, k_y)$    and  $g$ the function defined by:
\begin{equation}
\label{ggg}
2k_x g=i\ep k_y (k_x\Dy -k_y \Dx)u^{in}+2k_x u^{in}.
\end{equation}
The problem  may read as:
\begin{eqnarray}
\label{S-demi-2}
i(k_x \Dx + k_y \Dy)u + \frac{\ep}{ 2}(k_x^2 \Dyy -2k_x k_y
\Dxy + k_y^2 \Dxx)u+i\nu u=0, \qquad {\rm on } \; \cal{D},
\\
\label{CE-demi-2}
i\ep k_y (k_x\Dy -k_y \Dx)u+2k_x  u = 2 k_x g, \qquad \; {\rm on }
\Gamma_0.
\end{eqnarray}
In the sequel, the Fourier variables related to $x$ and $y$ respectively are $\xi$ and $\eta.$  The Fourier transform
 in $x$ and $y$ are denoted by ${\cal F}_x (\bullet)$ and ${\cal F}_y (\bullet)$, moreover ${\cal F}_y (u;x,.)$ denotes the Fourier transform of $u(x,.)$.

Here and in the sequel, $\sqrt{~~}$ denotes the principal determination of the square root (its real part is positive).
Denote:
$$
 R_- (i\eta)=i\frac{ k_x  \eta }{ k_y}-i\frac{ k_x }{ {\ep  k_y^2}}( 1 - \sqrt{1-2\frac{\ep k_y  \eta }{ k_x^2} + 2i \nu \frac{\ep k_y^2 }{ k_x^2}}).
$$
Since $\nu >0,$ one can define $R_-$ without ambiguity and
one checks that ${\cal R}e(R_- (i\eta) ) <0 $ for all $\eta.$
Let ${\cal S}'(\R)$ be the space of tempered distributions.

\begin{proposition} \label{demi-theo-2}

Assume that $g \in {\cal S}'(\R),$ then there exists a unique  distribution
 $u(x,.)$  continuous from $\R ^+$ into ${\cal S}'_y(\R),$ solution to Problem~(\ref{S-demi-2})(\ref{CE-demi-2}). It is given by:
\begin{equation}
\label{demi-exact}
{\cal F}_y (u;x,\eta)= \frac{2{\cal F}_y (g;\eta) }{ 1+ \sqrt{1-2\frac{\ep k_y \eta}{ k_x^2} + 2i\nu \frac{\ep k_y^2}{ k_x^2}}}e^{R_- (i\eta) x}.
\end{equation}
It  satisfies also:
$$\biggl(\Dx - R_- (i\eta)\biggr) {\cal F}_y (u;x,\eta)=0.$$

                                \end{proposition}

\begin{proof}

The principle is to take the Fourier transform in $y$ of the problem, and afterwards we shall consider Fourier transform in $x$ of the equation extended to the whole space.

Let $u$ be a solution of Problem~(\ref{S-demi-2})(\ref{CE-demi-2}) and $v$ the extension of $u$ by zero in the whole space: $v(x,y)= u(x,y) \3 $.
 By introducing formally the function  $v$ in Equation~(\ref{S-demi-2}) we get:
$$i\f{k}\cdot\f{\mathbf{\nabla }}v+\frac{\epsilon }{2}%
\Delta _{\bot }v+i\nu v= \biggl(\bigl(ik_x -\frac{ \ep k_y }{ 2}(2k_x \Dy -  k_y \Dx)\bigr)u(0,y)\biggl){\delta}_{x=0} + \frac{\ep k_y^2 }{ 2}u(0,y)
{\delta}^{'}_{x=0}.
$$
The term $\Dx u(0,y)$ is defined by the entrance boundary condition
(\ref{CE-demi-2}), so we get:
$$i\f{k}\cdot\f{\mathbf{\nabla }}v+\frac{\epsilon }{2}%
\Delta _{\bot }v+i\nu v=  ik_x g(y)\delta_{x=0} -\frac{\ep k_y }{ 2}\bigl( k_x \Dy u(0,y){\delta}_{x=0} - k_y u(0,y){\delta}^{'}_{x=0}\bigr).$$
Assuming that $u\in {\cal C}(\R_+, {\cal S}'(\R))$, we are allowed to take the Fourier transform of this expression. Let us define $P(X,Y)$ as
 the polynomial which characterizes the differential operator
 of the equation, that is to say:
 $$P(\Dx, \Dy)= i(k_x \dx + k_y \dy) + \frac{\epsilon }{2}(k_y^2 \Dxx-2k_x k_y \Dxy +k_x^2 \Dyy)+i\nu.  $$
Writing $u_0(y)=u(0,y)$, the Fourier transform in $y$ of the equation in $v$ reads:
$$P(\Dx, i\eta){\cal F}_y (v;x,\eta) =
 \frac{\ep k_y^2}{ 2}\biggl\{\biggl(\frac{2i k_x}{ \ep k_y^2}{\cal F}_y(g;\eta) - i \frac{k_x}{ k_y} \eta {\cal F}_y( u_0;\eta)\biggr)\delta_{x=0} + {\cal F}_y(u_0;\eta)\delta^{'}_{x=0}\biggr\}.$$
Polynomial $P$ may be factorized as:
\begin{equation}
P(\Dx,i\eta)=\frac{\ep k_y^2}{ 2}\biggl(\Dx -R_+ (i\eta)\biggr)\biggl(\Dx - R_- (i\eta)\biggr),
\label{fac1}
\end{equation}
where we define $R_\pm (i\eta)=i\frac{
 k_x}{ k_y}\eta -i\frac{ k_x}{ {\ep  k_y^2}}\biggl( 1 \pm \sqrt{1-2\frac{\ep k_y  \eta }{ k_x^2} + 2i \nu  \frac{\ep k_y^2}{ k_x^2}}\biggr).
$ Thus:
$$
\biggl(\Dx -R_+(i\eta)\biggr)\biggl(\Dx - R_-(i\eta)\biggr){\cal F}_y (v;x,\eta)=
$$\begin{equation}\label{new}
 \biggl(\frac{2i k_x }{ \ep k_y^2}{\cal F}_y(g;\eta) - i \frac{k_x}{ k_y} \eta {\cal F}_y( u_0;\eta)\biggr)\delta_{x=0} + {\cal F}_y(u_0;\eta)\delta^{'}_{x=0}.
\end{equation}
We now show that there is a unique acceptable solution for this ordinary differential equation.
Let us take its Fourier transform in $x$:
$$\biggl(i\xi -R_+(i\eta)\biggr)\biggl(i\xi - R_-(i\eta)\biggr){\cal F}_x {\cal F}_y (v;\xi,\eta)= \frac{2i k_x}{ \ep k_y^2}{\cal F}_y(g;\eta) - i (\frac{k_x}{ k_y} \eta -\xi) {\cal F}_y( u_0;\eta).$$
Since ${\cal R}e \bigl(i\xi - R_\pm (i\eta)\bigr) \neq 0$, we can divide each side of this equation 
by $\frac{2}{ \ep k_y^2} P(i \xi,i\eta):$ 

$${\cal F}_x {\cal F}_y (v;\xi,\eta) = \frac{\alpha^+ (\eta)}{ i\xi- R_+(i\eta)} + \frac{\alpha^-(\eta) }{ i\xi - R_-(i\eta)},$$
where
$\alpha^\pm(\eta)=\pm \frac{R_-(i\eta) -i\frac{k_x}{ k_y}\eta}{ R_+(i\eta) - R_-(i\eta)}{\cal F}_y(u_0;\eta) \pm \frac{2i k_x}{ \ep k_y^2} \frac1{ R_+(i\eta) - R_-(i\eta)} {\cal F}_y(g;\eta) .$



If $\theta \in \C \backslash \R$, one knows that:
$$\frac{1}{ i\xi-\theta} = \left\{\begin{array}{ll}
 {\cal F}_x (\3 e^{\theta x};\xi) & \mbox{ if ${\cal R}e(\theta )<0$}\\
 -{\cal F}_x (\5 e^{\theta x};\xi) & \mbox{ if ${\cal R}e(\theta )>0.$}
 \end{array}\right.$$
Here ${\cal R}e (R_+)=-{\cal R}e (R_-) >0.$ According to the previous
remark, since $v(x,.)=0$ for $x$ negative, one gets 
 $\alpha^+ (\eta)=0$    and 
 $$
 {\cal F}_y(u;x,\eta) = \alpha^- (\eta)e^{R_-(i\eta)x}\3,
 $$
so we get ${\cal F}_y(u_0;\eta) = -\frac{2i k_x }{ \ep k_y^2}\frac{{\cal F}_y(g;\eta)}{ R_+(i\eta) -i\frac{k_x}{ k_y}\eta}.$
Equality~(\ref{demi-exact}) and the last assertion follow.
\end{proof}

Notice that we can easily calculate, with this formula, the value of the derivative $\f{k}\cdot \nabla u$. As soon as $u$ is regular enough, we can perform an
asymptotic expansion in $\ep$ and $\nu$, and find: $\f{k}\cdot \nabla u= O(\ep +\nu)$.

From this result, one deduces the following stability  result.
\begin{coroll}\label{coroll-1}
If $g \in H^{-\frac{1}{2}} (\R)$ then the solution $u$ to Problem~(\ref{S-demi-2})(\ref{CE-demi-2}) is continuous from 
$\R ^+$ into $ L^2_y (\R)$,
 and it satisfies, for some constant~$C$ not depending on the coefficient $\nu$:
$$||u||_{L^{\infty}_x(\R_+, L^2_y (\R))} \leq C ||g||_{ H^{-\frac{1}{2} } (\R)}.$$
\end{coroll}

Since $C$ does not depend on the absorption coefficient $\nu,$ one can check that if $u^{in}$  is smooth enough,  for
 $x$ fixed,  the function $u(x,.)$ converges strongly to a function in
 $L^2_y$  when $\nu\to 0.$ Therefore, one may claim that there exists a bounded solution $u$ to Problem~(\ref{S-demi-2})(\ref{CE-demi-2}), even if $\nu=0$.

                \begin{proof}

Let us integrate with respect to $\eta$ the square modulus of both sides of Equation~(\ref{demi-exact}). Since $|e^{R_-(i\eta)x}|= e^{{\cal R}e(R_-(i\eta))x} \leq 1$ and:
$$ \int {|{\cal F}_y(g;\eta)|^2
   (1+|\eta|^2)^{-\frac{1}{2}}d\eta }= ||g||_{H^{-\frac{1}{2} }(\R)}^2 ,$$
it suffices to show
 that there exists a constant $C _1>0,$ not depending on $\nu$,
 such that:
\begin{equation}
1 + |\eta|^2  \leq C_1 
\bigg|1 +\sqrt{1 - \frac{2\ep k_y}{ k_x^2} \eta
+ 2i\ep \nu \frac{k_y^2}{ k_x^2}} \bigg|^4 \qquad \forall \eta\;\in \R.
\label{lemde}
\end{equation}
So, if we denote $X= 1 - \frac{2\ep k_y}{ k_x^2} \eta$ and $N=2\ep \nu \frac{k_y^2}{ k_x^2}$, one first sees that:
$$|1 + \sqrt{X +iN}|^2 = 1 + \sqrt{X^2 + N^2} +2(X^2 + N^2)^{\frac{1}{4}} cos(\frac{\pi}{4} - \frac{{\rm Argtan}{X/N}}{2}) \geq \sqrt{1 + X^2}$$
(indeed the cosine is nonnegative).
With $a=\frac{k_x^2}{2\ep k_y},$ we have $1+|\eta|^2=1+a^2(1-X)^2$ and it is easy to check that  $ 1+a^2(1-X)^2 \leq C_1 (1+ X^2)$ for $C_1=2a^2 +1$~; Inequality~(\ref{lemde}) follows.
\end{proof}

\

{\bf Remark:} with the same techniques, one can also find existence and uniqueness of a solution in other spaces, for instance, if $ \frac{{\cal F}_y (g;\eta) }{ (1+ |\eta|^2)^{1/8}  }  \in L^2_\eta(\R),$ we have  $u \in L^2 (\mathcal{D}).$

Since $| {\cal F}_y (g;\eta) | \le  C (1+
|\eta|^2)^{1/2}
|{\cal F}_y (u^{in};\eta)|,$ that means that
if $u^{in}$  is smooth enough (in $H^{3/4}$ for example), the solution
$u$  belongs to $L^2(\mathcal{D} ).$

  \subsection { Remark on the Problem on the Quadrant}
\label{subsection-quart}

We now consider the same problem (\ref{S-demi-2})(\ref{CE-demi-2})
but restricted to the quadrant $\{(x,y) \; {\rm s.t.}\; 
x\ge 0,\; y\geq 0 \}$. To find a
good  absorbing boundary condition on the boundary $\{y=0\}$, 
we formally factorize the differential operator of Equation~(\ref{S-demi-2}) as follows: 
\begin{equation}  
P(\Dx,\Dy)= \ep \frac{k_x^2}{2}\bigl(\Dy -A_+(\Dx)\bigr)\bigl(\Dy -A_-(\Dx)\bigr),
\label{fact}
\end{equation} 
where $A_+(.)$ and $A_-(.)$ are the roots of $P$ considered as 
polynomials in $\Dy:$ 
$$
A_\pm (\Dx)= \frac{k_y}{k_x} \Dx -i \frac{k_y}{ \ep k_x^2}\bigl(1
\pm \sqrt{1+\frac{2i\ep k_x}{ k_y^2} \Dx +2i\ep\nu \frac{k_x^2}{
k_y^2}}\bigr)=
\frac{k_y}{ k_x} \Dx -i \frac{k_y}{ \ep k_x^2}
\mp \frac{1}{ \ep k_x^2}  \sqrt{-k_y^2-2i\ep k_x \Dx -2i\ep\nu k_x^2 }.
$$
The definition of the fractional derivative is classical and  is based on Fourier transform.
The quadrant problem that we consider consists of Equations~(\ref{S-demi-2})(\ref{CE-demi-2}) supplemented with the following boundary condition
\begin{equation}\label{CT}
\Dy u - A_+ (\Dx) (u) = 0, \quad  \forall x>0, \;
{\rm for} \quad  y=0.
\end{equation}


Then, we have the following result, which is detailed in \cite{D,AUTRE} (for related boundary value problems for classical Schr\"odinger equations, see for example \cite{Arnold}).
           \begin{proposition} \label{quart-theo1}
Assume $g \in H^{-\frac{1}{2}}(\R^+)$ and its support is in $ (0,+\infty).$
 Let $u$ be the solution of the
 half-space problem~(\ref{S-demi-2})(\ref{CE-demi-2}). There is a unique
 solution $U$ continuous from 
$\R ^+$ into $ L^{2}_y (\R^+)$ of Problem~(\ref{S-demi-2})(\ref{CE-demi-2})(\ref{CT})  and it satisfies

i) if $k_y>0$, then $U=u \4$, 

ii) if $k_y <0$ and if the incoming data is given by $g(y)=h(y-a)$ with $a>0$, then:

$$
\lim\limits_{a\to +\infty} ||U- u\4||_{L^\infty(\R^+,L^2_y(\R^+))}=0.
$$
\end{proposition}

\section{Numerical Scheme}

\label{section-num}


Let us consider the domain: 
\begin{equation*}
\mathcal{D}=\{(x,y) : \; 0\leq x \leq L_x, \; y_0\leq y \leq y_0+L_y\}.
\end{equation*}
On this domain, we address the numerical solution of the following equation: 
\begin{equation}\label{eq:advSchro:num}
i(k_x {\partial_x} + k_y {\partial_y})u + {\frac{\epsilon }{2}}\Delta_{\bot}
u+i\nu u -\mu u =0,
\end{equation}
where $\nu = \nu (\mathbf{x})$ and $\mu = \mu (\mathbf{x});$ it is
supplemented by the same boundary condition as before on $\{x=0\}$ : 
\begin{equation*}
i\epsilon k_y (k_x{\partial_y} -k_y {\partial_x})u +2k_x u =2 k_x g,
\end{equation*}
where $g$ is given by Equation~(\ref{ggg}). It is the same problem as in Section \ref%
{section-demi-2}, except that the coefficients $\nu$ and $\mu$ may be
functions of $\mathbf{x}.$ In the sequel, we consider alternatively the case where 
$\mu$ is a function of $|u|;$ as a matter of fact, we can take 

\begin{equation*}
\mu= f(|u|),\qquad \mathrm{where} \; f(w)=e^{-\alpha w^2}-1,
\end{equation*}
with $\alpha $ a positive constant (for a justification of this model, see for example \cite{rose} \cite{sen} ).

The interesting problems involve a very small coefficient $\nu$, and it may
be necessary to have $\alpha$ sufficiently small so that there is no blow-up
of the solution.

\subsection{Description of the Scheme} \

\label{subsection-schema}

Let us set : 
\begin{equation*}
\nu=\nu_0 + \nu_1 \; \; \; \; \mathrm{with }\; \; \; \; \nu_0= \inf \nu,
\end{equation*}
so $\nu_0$ is a constant and $\nu_1$ a function of $\mathbf{x}.$
One discretizes the problem according to a regular grid, we denote by $%
\delta x ,\; \delta y$ the space step in the two directions and by $n$ and $%
j $ the indices corresponding respectively to $x$ and $y;$ 
then $u^n_j\approx u(n\delta x,j\delta y).$

The numerical method is based on a space marching technique according to the 
$x$ variable and a splitting with respect to this variable. According to
Proposition \ref{demi-theo-2}, when the value of $u^n$ is known, it would be
possible to evaluate a first intermediate value $u^{\mathrm{inter}}$ by
solving on $[x^n, x^{n }+\delta x ] $ the following equation: 
\begin{equation*}
( k_x {\partial_x} +k_y {\partial_y}) u -i \frac{\epsilon}{2} \Delta_{\perp}
u+ \nu_0 u = 0.
\end{equation*}
it would be given by $\mathcal{F} ( u^{\mathrm{inter}} )=\mathcal{F} (u^n)
e^{R_- (i\eta) \delta x}$ (here we denote ${\cal F}={\cal F}_y$).

As a matter of fact, in order to have an accurate treatment of the advection
term, we prefer to perform the following simple splitting : at each space
step $[x^n, x^{n }+\delta x ],$ one solves succesively
\begin{eqnarray*}
k_x {\partial_x} u -i \frac{\epsilon}{2} \Delta_{\perp} u+ \nu_0 u = 0, \\
k_x {\partial_x} u + k_y {\partial_y} u + (\nu_1 +i \mu ) u = 0.
\end{eqnarray*}

\subsubsection{Initialization}

For the initial condition, recall that
\begin{equation*}
g=i\epsilon {\frac{k_{y}}{2k_{x}}}(k_{x}{\partial _{y}}-k_{y}{\partial _{x}}%
)u^{in}+u^{in},
\end{equation*}
where the input data $u^{in}=u_{|x=0}^{in}$ is a smooth function of the
transverse variable $Y=\mathbf{k}_{\bot }\cdot \mathbf{x}=k_{x}y-k_{y}x$
which values zero around the corner points $y=y_{0}$ and $y=y_{0}+L_{y},$ so
one can take its Fourier transform.

To determine the boundary value $u^{0}$ of $u,$ we use Formula~(\ref%
{demi-exact}) 
\begin{equation}
\mathcal{F}(u^{0})={\frac{2\mathcal{F}(g)}{1+\sqrt{1-2{\frac{\epsilon
k_{y}\eta }{k_{x}^{2}}}+2i\nu ^{in}{\frac{\epsilon k_{y}^{2}}{k_{x}^{2}}}}}}.
\label{init-u0-g}
\end{equation}

That is to say, $(u_{j}^{0})_{j}$ is obtained by taking the FFT (Fast
Fourier Transform) of $g$, dividing this function of $\eta $ by the function 
$1+\sqrt{1-2{\frac{\epsilon k_{y}\eta }{k_{x}^{2}}}+2i\nu ^{in}{\frac{%
\epsilon k_{y}^{2}}{k_{x}^{2}}}}$ and then taking the IFFT (Inverse Fast
Fourier Transform) of the result.

\bigskip Generally, the input data $u^{in}$ is a sum of Gaussian functions
whose half-height width is in the order of a characteristic length $L_{s}$
which is the typical value of the speckle width (a speckle is a hot spot
inside the laser beam) and $L_{s}$ is generally larger than $20$ times $%
\varepsilon .$ Then one checks that for values of $\epsilon /L_{s}$ less
than $0.1,$ the term $i\epsilon k_{y}(k_{x}{\partial _{y}}-k_{y}{\partial
_{x}})u^{in}$ that appears in the previous formula for $g$ is a corrective
term and it is possible to take simply $g$ equal to $u^{in}.$

\subsubsection{First stage: Fourier transform}

\label{stages1-2}

The first stage is to solve 
\begin{equation}
k_{x}{\partial _{x}}u-i\frac{\epsilon }{2}\Delta _{\perp }u+\nu _{0}u=0,
\label{SNaa}
\end{equation}%
and we proceed from $u^{n}$ to $u^{n\#}$. Practically,  from Proposition \ref{demi-theo-2},  we get immediately : 
\begin{equation*}
\mathcal{F}(u^{n\#})=\mathcal{F}(u^{n})e^{(R_{-}(i\eta )+i\eta {\frac{k_{y}}{%
k_{x}}})\delta x}.
\end{equation*}%
In fact, we have 
\begin{equation}
R_{-}(i\eta )+i\eta {\frac{k_{y}}{k_{x}}}=-{\frac{2\nu _{0}}{k_{x}(1+\sqrt{%
1-2{\frac{\epsilon k_{y}\eta }{k_{x}^{2}}}+2i\nu _0 {\frac{\epsilon k_{y}^{2}}{%
k_{x}^{2}}}})}}-{\frac{2i\eta \epsilon (\eta -i\nu _{0}k_{y})}{k_{x}^{3}(1+%
\sqrt{1-2{\frac{\epsilon k_{y}\eta }{k_{x}^{2}}}+2i\nu _{0}{\frac{\epsilon
k_{y}^{2}}{k_{x}^{2}}}})^{2}}}.  \label{R-epky}
\end{equation}%
Notice that this formula may be used even if $\nu _{0}$ is equal to zero, provided that the square root of the complex quantity is well defined.

So, after a FFT on $(u_{{}}^{n})$, we multiply it by $e^{(R_{-}(i\eta
)+i\eta {\frac{k_{y}}{k_{x}}})\delta x}$ and then apply an inverse FFT. We
denote $\bigl(u_{j}^{n\#}\bigr)$ the value of the intermediate function, in
the cell $(n,j)$.

\subsubsection{Second stage: finite difference scheme}
\label{subsubsec:scheme:second}
\textbf{ Boundary conditions on the edges }$\{y=0\}$ \textbf{and} $%
\{y=L\}$

It is well known that for this kind of propagation model, the boundary
treatment is sensitive; see for example \cite{ber} for the case of wave
equations. In our case the problem is somehow different since there is a
privileged direction of propagation: as we use a FFT technique, the key
point at each stage of the space marching scheme is to force the values of
the numerical solution to be negligeable on both edges. Therefore we use a
damping method which is well known by physicists who address this kind of
problem \cite{had}. The principle is to introduce in a strip near each edge an
artificial absorbing coefficient denoted by $B$; it decreases progressively
on the first five cells near the edge and is very large on the edge. More precisely, if $\nu _{1,j}^{n}$ denotes the
value of $\nu _{1}$ in cell $(n,j),$ one replaces $\nu _{1,j}^{n}$ by $\nu
_{1,j}^{n}+B_{j}$ where the artificial coefficient $B_{j}$ is defined by

\begin{equation}
\begin{array}{llll}
B_{j} &=&b \beta ^{5-j}& \text{if }j\leq 5 \\
&=&b \beta ^{5-J_{\max }+j}& \text{if }J_{\max }-j\leq 5 \\
&=&0 \qquad \text{\ }& \text{elsewhere,}
\end{array}\label{eq:defB}
\end{equation}
with $\beta $  typically in the order of $10$ to $100.$ The numerical tests below (with a
characteristic value of $b$ in the order of $0.1$ to $1$) show that this
technique leads to get a vanishing value of the solution on the edges.
One checks on Table \ref{tab:B} that the value of the solution (outside the artificial absorbing layers) is almost independant from the choosen values of $b$ and $\beta.$ Indeed, near the boundary, the main step is the advection one and it is crucial to have a numerical solution which is negligible near the boundary cell, in order to avoid a spurious ray to appear on the opposite boundary, due to the FFT.
Notice that, according to the advection scheme by space marching, the modification in the artificial layer at position $x^n$ has no significant impact on the value outside the artificial layer at position $x^{n+1} .$

\textbf{First order scheme.}

In this stage, we solve on $[x^{n},x^{n}+\delta x]$ the following equation: 
\begin{equation}
k_{x}{\partial _{x}}u+k_{y}{\partial _{y}}u+\nu _{1}(x^{n})u+i\mu u=0.
\label{SNbb}
\end{equation}%
 To do this, we use
standard finite difference methods. Assume that $k_{y}>0$ (the case $k_{y}<0$
is similar). We consider an  upwind method, given that the CFL
stability criteria $\theta \leq 1$ must be checked, where 
$$\theta ={\frac{k_{y}}{k_{x}}}{\frac{\delta x}{\delta y}}. 
$$

 The initial value is now $u_{j}^{n\#}$ and we get the final value $u_{j}^{n+1}$
 by setting 
\begin{equation}
{\frac{k_{x}}{\delta x}}(u_{j}^{n+1}-u_{j}^{n\#})+{\frac{k_{y}}{\delta y}}
(u_{j}^{n\#}-u_{j-1}^{n\#})+\biggl(\nu _{1,j}^{n}+i\mu_{j}^{n}   \biggr)
\bigl( \frac{u_{\theta j}^{n\#}+u_{j}^{n+1}} {2} \bigr)+B_{j}u_{j}^{n+1}=0,
\label{scheme1}
\end{equation}
where $u^{n\#}_{\theta j}= \theta u^{n\#}_{j-1} + (1-\theta)u^{n\#}_j. $ It is the value of the function on the characteristic line passing by $(x^{n+1},y_j);$ 
for the first cell, we set  $u^{n\#}_{-1}=0$. 

For the nonlinear model where the term  $\mu $ is replaced by  $f(|u|)$, the coefficient $\mu_{j}^{n}$  has to be
replaced by $f(|u^{n\#}_{\theta j}|)$ .  

\ 

\textbf{Second order scheme}

\ 

When $\theta =1$, the previous scheme gives very accurate results, but in
real cases it is not possible to impose this condition, one has $\theta <1$
and results are much worse (see Table \ref{tab:conv:ord1:CFL}). We improve the numerical
scheme when $\theta <1$ by using a second order scheme as in all advection
problems. To do this, we choose a flux-limiter method (see \cite{RL}), with
the Van Leer function as limiter (tests prove it to be the best one: see Figure \ref{fig:Ord2Flux} and Section \ref{subsubsec:conv}). That
is to say, we introduce the function $\phi $ which depends on the ratio $%
\lambda $ of the gradient of the function $u^{\#}$ in two neighboring
cells: 
\begin{equation}\label{eq:phi:VanLeer}
\phi (\lambda )={\frac{|\lambda |+\lambda }{1+|\lambda |}}.
\end{equation}%
We have to solve simultaneously two scalar equations (one for the real and
one for the imaginary part) with the same flux limiter, so we have to choose
one single significant quantity to estimate the flux limitor: we choose the
energy of the laser, \emph{i.e.} $|u|^{2}$, and evaluate $\phi $ in terms of 
$|u_{j}|^{2}$ and not of $|u_{j}|$: 
\begin{equation*}
\lambda _{j}={\frac{|u_{j}^{\#}|^{2}-|u_{j-1}^{\#}|^{2}}{|u_{j+1}^{%
\#}|^{2}-|u_{j}^{\#}|^{2}}}.
\end{equation*}%
We now replace, in the first order scheme, the term derivative in $y$, $%
u_{j}^{\#}-u_{j-1}^{\#},$ by $F_{j}-F_{j-1}$ where the flux $F_{j}$ is
defined as: 
\begin{equation*}
F_{j}=u_{j}^{\#}+{\frac{1}{2}}(1-\theta )(u_{j+1}^{\#}-u_{j}^{\#})\phi
(\lambda _{j}).
\end{equation*}%
The second order scheme is now: 
\begin{equation}
{\frac{k_{x}}{\delta x}}(u_{j}^{n+1}-u_{j}^{n\#})+{\frac{k_{y}}{\delta y}}%
(F_{j}^{n}-F_{j-1}^{n})+\biggl(\nu _{1,j}^{n}+i \mu_{j}^{n} \biggr)\bigl({%
\frac{u_{\theta j}^{n\#}+u_{j}^{n+1}}{2}}\bigr)+B_{j}^{{}}u_{j}^{n+1}=0.
\label{scheme2}
\end{equation}

\subsubsection{Numerical method for two-ray model} 
 \label{subsubsec:2ray}
One may also consider a more complex model with two rays crossing each
other, with two different propagation vectors $\mathbf{k}^1$ and $\mathbf{k}%
^2$ (one with positive and one with negative $y-$component: $k_y^1>0$ and $%
k_y^2<0$.) To do so, it is necessary to evaluate the nonlinear term $f(|u|).$
Theoretically, the laser energy is: 
\begin{equation*}
|\Psi|^2=|u^1e^{i{\frac{\mathbf{k}^1\cdot \mathbf{x}}{\epsilon}}}+u^2e^{i{%
\frac{\mathbf{k}^2\cdot \mathbf{x}}{\epsilon}}}|= |u^1|^2+|u^2|^2+ 2 
\mathcal{R}e\bigl( u^1u^{2*} e^{i{\frac{(\mathbf{k}^1-\mathbf{k}^2)}{\epsilon%
}}\cdot \mathbf{x}}\bigr).
\end{equation*}
But we are in the framework of W.K.B. approximation and we do not model the
fluctuation of the solution at the wavelength level. Hence, the term $f$ has
to be taken on a function $w$ corresponding to the variation of the index of
refraction, which is here the average value of $|u|$ over a wavelength: 
\begin{equation*}
w=\sqrt{|u^1|^2+|u^2|^2}.
\end{equation*}
One considers the following model, for $p=1,2$: 
\begin{equation*}
i\mathbf{k}^p\cdot\nabla u^p + {\frac{\epsilon}{2}} \Delta^p_{\bot}+i\nu u^p
= f(\sqrt{|u^1|^2+|u^2|^2})u^p.
\end{equation*}
The first stage of the previous scheme is the same as before : for each ray,
we consider Equation (\ref{SNaa}) with its own propagation direction $%
\mathbf{k}^1$ or $\mathbf{k}^2$. The interaction between the two rays
changes only the nonlinear term of the second stage.

\subsection{Properties of the scheme}
\label{subsec:scheme:prop}
\subsubsection{Stability}\label{subsubsec:stability}

Let us denote $||v^n||^2_{l^2}=\sum\limits_j|v^n_j|^2 \delta y.$
\begin{proposition}\label{th:stability}
The numerical first order scheme is monotone decreasing for the $l^2$-norm, 
\emph{i.e.} the following inequality stands
\begin{equation}
\forall\, n\in {\mathbf{N}}\,,\qquad ||u^n||_{l^2} \leq
||u^{n+1}||_{l^2}.
\end{equation}
Moreover, the previous inequality is strict if $\nu \neq 0.$
\end{proposition}

\begin{proof}
\begin{enumerate}
\item First stage: Fast Fourier Transform

Let us denote by $\zeta$ the discrete variable associated to $\eta$. On the
one hand, since 
\begin{equation*}
u^{n\#}=IFFT \biggl(e^{\bigl(R_-(i\zeta)+i\zeta{\frac{k_y}{k_x}}\bigr)\delta
x} FFT(u^n) \biggr)
\end{equation*}
and since the FFT conserves the $l^2$-norm, we have: 
\begin{equation*}
||u^{n\#}||_{l^2}=||e^{\bigl(R_-(i\zeta)+i\zeta{\frac{k_y}{k_x}}\bigr)\delta
x} FFT(u^n)||_{l^2}.
\end{equation*}
On the second hand, the inequality $\mathcal{R}e \bigl(R_-(i\zeta)\bigr)\leq
0$ implies that 
\begin{equation*}
|e^{\bigl(R_-(i\zeta)+i\zeta {\frac{k_y}{k_x}}\bigr)\delta x}|\leq 1,
\end{equation*}
with an equality \emph{iff} $\nu_0=0.$ We deduce that: 
\begin{equation*}
||e^{\bigl(R_-(i\zeta)+i\zeta{\frac{k_y}{k_x}}\bigr)\delta
x}U(\zeta)||_{l^2}\leq ||U(\zeta)||_{l^2},
\end{equation*}
and conclude: 
\begin{equation*}
||u^{n\#}||_{l^2} \leq ||u^n||_{l^2},
\end{equation*}
with $||u^{n\#}||_{l^2} = ||u^n||_{l^2}$ \emph{iff} $\nu_0=0.$

\item Second stage: upwind scheme

For the first order scheme, Relation~(\ref{scheme1}) gives us that: 
\begin{equation*}
u_{j}^{n+1}={\frac{{\frac{k_{x}}{\delta _{x}}}u_{j}^{n\#}-{\frac{k_{y}}{%
\delta _{y}}}(u_{j}^{n\#}-u_{j-1}^{n\#})-{\frac{1}{2}}\bigl(\nu
_{1,j}^{n}+i \mu_{j}^{n}  \bigr)u_{\theta _{j}}^{n\#}}{{\frac{k_{x}}{\delta
_{x}}}+{\frac{1}{2}}\bigl(\nu _{1,j}^{n}+i  \mu_{j}^{n} \bigr)+B_{j}}}
\end{equation*}%
Provided that ${\frac{k_{x}}{\delta _{x}}}u_{j}^{n\#}-{\frac{k_{y}}{\delta
_{y}}}(u_{j}^{n\#}-u_{j-1}^{n\#})={\frac{k_{x}}{\delta _{x}}}u_{\theta
_{j}}^{n\#},$ we obtain: 
\begin{equation}\label{eq:step2}
u_{j}^{n+1}={\frac{{\frac{k_{x}}{\delta _{x}}}-{\frac{1}{2}}\bigl(\nu
_{1,j}^{n}+i \mu_{j}^{n} \bigr)}{{\frac{k_{x}}{\delta _{x}}}+{\frac{1}{2}}%
\bigl(\nu _{1,j}^{n}+i  \mu_{j}^{n}  \bigr)  +B_{j}}}u_{\theta _{j}}^{n\#}.
\end{equation}%
Since the modulus of the multiplicative coefficient in the right-hand side is smaller than one,
this leads to $||u^{n+1}||_{l^{2}}\leq ||\bigl(u_{\theta _{j}}^{n\#}\bigr)%
_{j}||_{l^{2}}.$ By the triangle inequality: 
\begin{equation*}
||\bigl(u_{\theta _{j}}^{n\#}\bigr)_{j}||_{l^{2}}\leq \theta ||\bigl(%
u_{j-1}^{n\#}\bigr)_{j}||_{l^{2}}+(1-\theta )||\bigl(u_{j}^{n\#}\bigr)%
_{j}||_{l^{2}}\leq ||u^{n\#}||_{l^{2}},
\end{equation*}

which concludes the proof.
\end{enumerate}
\end{proof}

In the linear case, that is the case where $\mu$ is a data and not a function of $|u|,$ the scheme is obviously consistent, so Proposition \ref{th:stability}
implies the convergence of the scheme.

Concerning the second order scheme modifying the advection step, it is well known (cf \cite{RL}) that the effect of this technique with a flux-limiter is to allow small $CFL-$numbers with a better accuracy (than the first order scheme)
 without generating spurius oscillations. These assertions will be confirmed by numerical tests we have performed (see Section \ref{subsubsec:conv}).  

\subsubsection{Comparison with the classical Schr\"odinger equation}

If $k_y \rightarrow 0,$ Equation (\ref{eq:advSchro:num})
reduces to the classical Schr\"odinger equation, in the case $\mu=f(|u|):$ 
\begin{equation}  \label{schroky0}
i\partial_x u +{\frac{\epsilon}{2}}\partial^2_{yy}u +i\nu u -f(|u|)u=0,
\end{equation}
with a very simple boundary condition (notice that $g \rightarrow u^{in}$) 
\begin{equation}
u_{|x=0}=u^{in}.  \label{schroky0-bord}
\end{equation}

\begin{proposition}
If $k_y\to0,$ the solution given by the numerical scheme converges to the
solution of the classical Schr\"odinger problem~(\ref{schroky0} )(\ref%
{schroky0-bord}).
\end{proposition}

\begin{proof}

* Initializing. Formula (\ref{init-u0-g}) used in the scheme shows that 
\begin{equation*}
\lim\limits_{k_y\to 0}\mathcal{F}(u;x=0)=\mathcal{F}(g),
\end{equation*}
so the boundary condition tends to $u_{|x=0}=g,$ which is Equation (\ref%
{schroky0-bord}).

* First stage.
If $k_y$ tends to zero, \emph{i.e} when the ray tends to be perpendicular to
the boundary, Formula (\ref{R-epky}) shows that: 
\begin{equation*}
\lim\limits_{k_y\to 0} R_-(i\eta)+i\eta{\frac{k_y}{k_x}} = -\nu-i{\frac{%
\epsilon}{2}}\eta^2,
\end{equation*}
so $u^{n\#}$ given by the first stage is the solution of the classical
Schr\"odinger equation without potential: 
\begin{equation*}
i\partial_x u +{\frac{\epsilon}{2}}\partial^2_{yy}u +i\nu u =0,
\end{equation*}
which is the limit of the advection-Schr\"odinger equation.

* Second stage.
It corresponds to a classical discretization of the ordinary differential
equation: 
\begin{equation*}
\partial_x u +\nu_1 u +if(|u|)u=0.
\end{equation*}
In other words, the scheme is a classical splitting between dispersion and
refraction in the Schr\"odinger equation (\ref{schroky0}).
\end{proof}

\subsection{Numerical results}

\label{subsection-num-result}

Let us recall that the laser energy density is equal to $|u|^{2}.$ Moreover, 
the physical meaning of the absorption coefficient $\nu$ is the following: with a constant value of $\nu ,$ if there would be no diffraction operator, the laser intensity (integrated on a line orthogonal to the propagation direction) would decrease by a factor $1/e^2$ on a propagation
distance equal to $1/ \nu.$

\begin{figure}
\begin{minipage}[b]{0.46\linewidth}
\centering\epsfig{figure=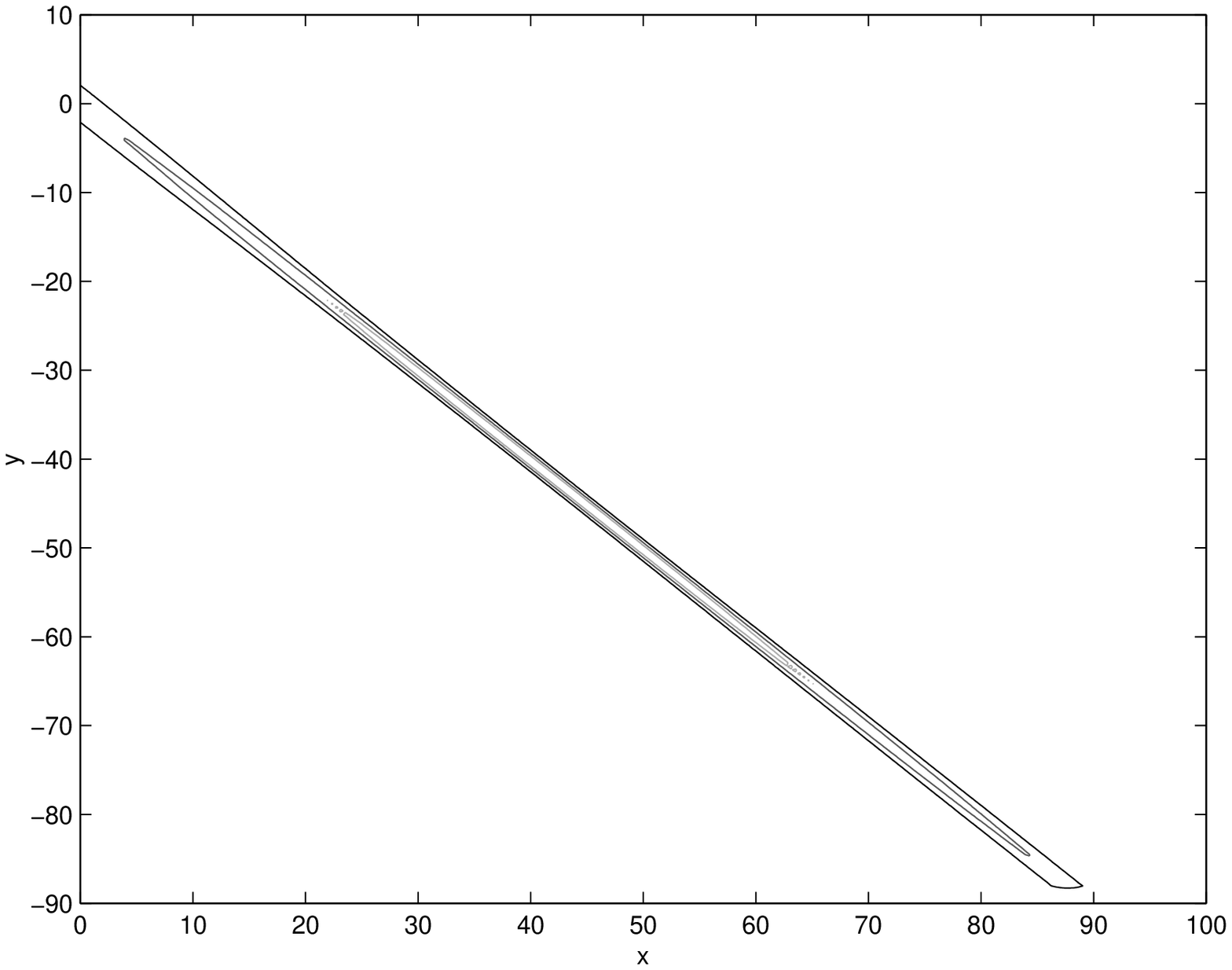, width=\linewidth, height=\linewidth}
\caption{\label{fig1-ref}Reference case: $\delta x=\delta y=0.05$, $CFL=1$. Then $L_{foc}=59.7,$ $Max(|u|^2)=2.14.$}
\end{minipage} \hfill
\begin{minipage}[b]{0.46\linewidth}
\centering\epsfig{figure= 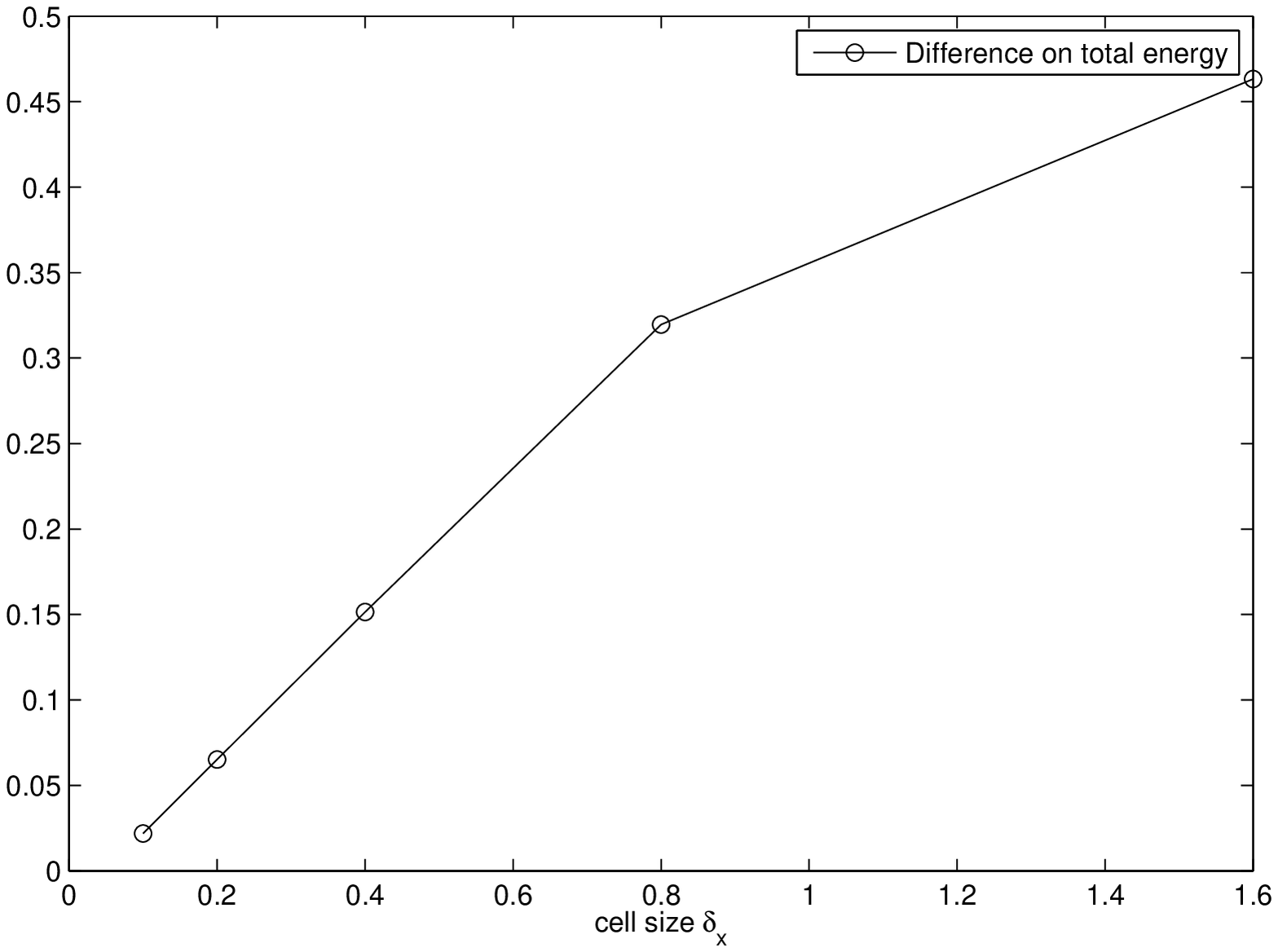, width=\linewidth, height=\linewidth}
\caption{\label{fig2-cvgce}1st order scheme convergence with $CFL=1$ as a function of cell size $\delta x$ (see Table \ref{tab:conv:ord1}).}
\end{minipage}
\end{figure}

We now give the standard numerical values used for the numerical tests.

\begin{enumerate}

\item
 For the incoming boundary  condition on the edge $x=0$, we take a Gaussian of amplitude $1$ centered at a point $(0,y_0)$ i.e. $%
u^{in}=\mathrm{exp}(-(k_x (y-y_0)-k_y x )^2 /{L_s^2})$ with $L_s=2.5 \;\mu m;$ which corresponds to the typical half-width of a speckle of a
 laser beam.

\item For the incidence angle, we take $-45^0,$ then  $\mathbf{k}%
=(-{\frac{\sqrt{2}}{2}},{\frac{\sqrt{2}}{2}}).$  

\item
$\epsilon=0.05 \; \mu m$, the wavelength of the laser is $2 \pi
\epsilon\approx 0.31 \; \mu m$.

\item
$\nu_0=\nu_1=5.10^{-4}  \; \mu m^{-1}.$  Notice that the larger the absorption coefficient, the easier the numerical simulation (indeed the laser energy decreases faster with respect to the propagation distance). 

\item
We take $\alpha=5.10^{-2}$. It depends on the electronic density of the
plasma: in the vacuum $\alpha$ would be null. This size order corresponds
either to a dense plasma or to a high laser intensity - since we have taken
a normalized value of the intensity corresponding to a maximum value of  $u^{in}$  equal to 1.

\item 
For the definition of the boundary layer $B,$ given by (\ref{eq:defB}), we take $b=0.1$ and $\beta=50.$
\end{enumerate}
All our figures represent the laser energy $|u|^{2}$.

 To be easier to read, our examples are variations with respect to the  case defined by the previous numerical values of the coefficients
 and  computed with a CFL number $\theta$ equal to 1 (see Figure \ref{fig1-ref}). With these assumptions, the scheme converges very well as the discretization step decreases 
 (see Table \ref{tab:conv:ord1}). 
 Due to the $\alpha$ coefficient, focusing occurs: the beam focuses and reaches a maximum, then decreases.  Notice that it may even focus several times for larger values of $\alpha.$
All our comparisons are made with this reference case, denoted $u^{\rm ref},$ in the fully converged situation (with mesh size $\delta x=0.05,$ corresponding to $2^{11}$ points on a domain length $L_x=100.$)

\subsubsection {Convergence of the scheme}
\label{subsubsec:conv}
{\bf Convergence of the first order scheme}

\begin{table}
\begin{center}
\begin{tabular}{|l|ccccc|c|}
\hline
\bf{Number of points} & $2^6$ & $2^7$  & $2^8$ & $2^9$ & $2^{10}$ & $\bf{2^{11}}$  \\ \hline
\bf{Mesh size $\delta x =\delta y$} & 1.6	& 0.8	& 0.4 &	0.2 &	0.1 &	\bf{0.05}  \\ \hline
\bf{Error on energy} $\Sigma_{j,n} ||u_j^n|^2-|u^{\rm ref,n}_j|^2|\delta x \delta y/|u^{\rm ref}|^{2} $ & 46 \% &	32\% &	15\% &	6\% &	2\%&	\bf{-}  \\ \hline
\bf{Focusing distance $L_{foc}$} & 82.7 & 61.4 &	59.5&	59.4	& 59.9 &	\bf{59.7}  \\ \hline
\bf{Error on focusing distance} & 38\% &	2.9\% &	0.4\% &	0.6\%&	0.3\%	& \bf{-}   \\ \hline
\bf{Maximum of energy $Max_{n,j}(|u^n_j|^2)$} & 1.74 &	2.16	& 2.13	& 2.13&	2.14&	\bf{2.14}  \\ \hline
\bf{Error on the maximum of energy} 
 & 19\% &	0.7\%	&0.4\%	& 0.4\%	& 0.07\%&	\bf{-}  \\ \hline
\end{tabular}

\caption{Convergence of the scheme, with $CFL=1.$ The last column represents the fully converged reference case $u^{\rm ref}.$
\label{tab:conv:ord1}}
\end{center}
\end{table}
\noindent
We first take the CFL number equal to 1, which is the case where the first and the second order schemes are equivalent. 
To verify the convergence of the scheme, we have three possible indicators. 
 A first indicator is the total energy in the physical domain of interest (that is to say, outside the artificial absorbing layer) which is equal to the $l^1-$norm of the energy: we denote it by
$$|u|^2 = \Sigma_{n,j} |u_j^n|^2 \delta x \delta y.$$
 So we compare this quantity to the corresponding one of the fully converged case $|u^{\rm ref}|^2 $; in the two first tables, we give the values of the relative error  $ \Sigma_{n,j} ||u_j^n|^2-|u^{\rm ref,n}_j|^2 |\delta x \delta y /|u^{\rm ref}|^2 $ for different cases.
 Now, if we want to compare for instance the effects of the variation of the incidence angle, two other indicators are more relevant in the framework of  the nonlinear model.
 One is given by the focusing distance: we can look for the focusing maximal point $L_{foc}$ and we measure the distance from $L_{foc}$ to the origin of the ray. A last indicator is the maximal value of the energy. 
 These last two indicators are quite sensitive.  For the nonlinear model,  the numerical results are  illustrated by Figure \ref{fig2-cvgce}  for the reference case ; the  estimates of the indicators are close to the ones of the reference case  when the spatial step decreases (see Table \ref{tab:conv:ord1}). 

Thus, we may conclude that when $CFL=1,$ we reach an accurate result even for $\delta x=\delta y=0.4,$ and that the focusing phenomenon is very well captured.

\begin{figure}
\begin{minipage}[b]{0.46\linewidth}
\centering\epsfig{figure=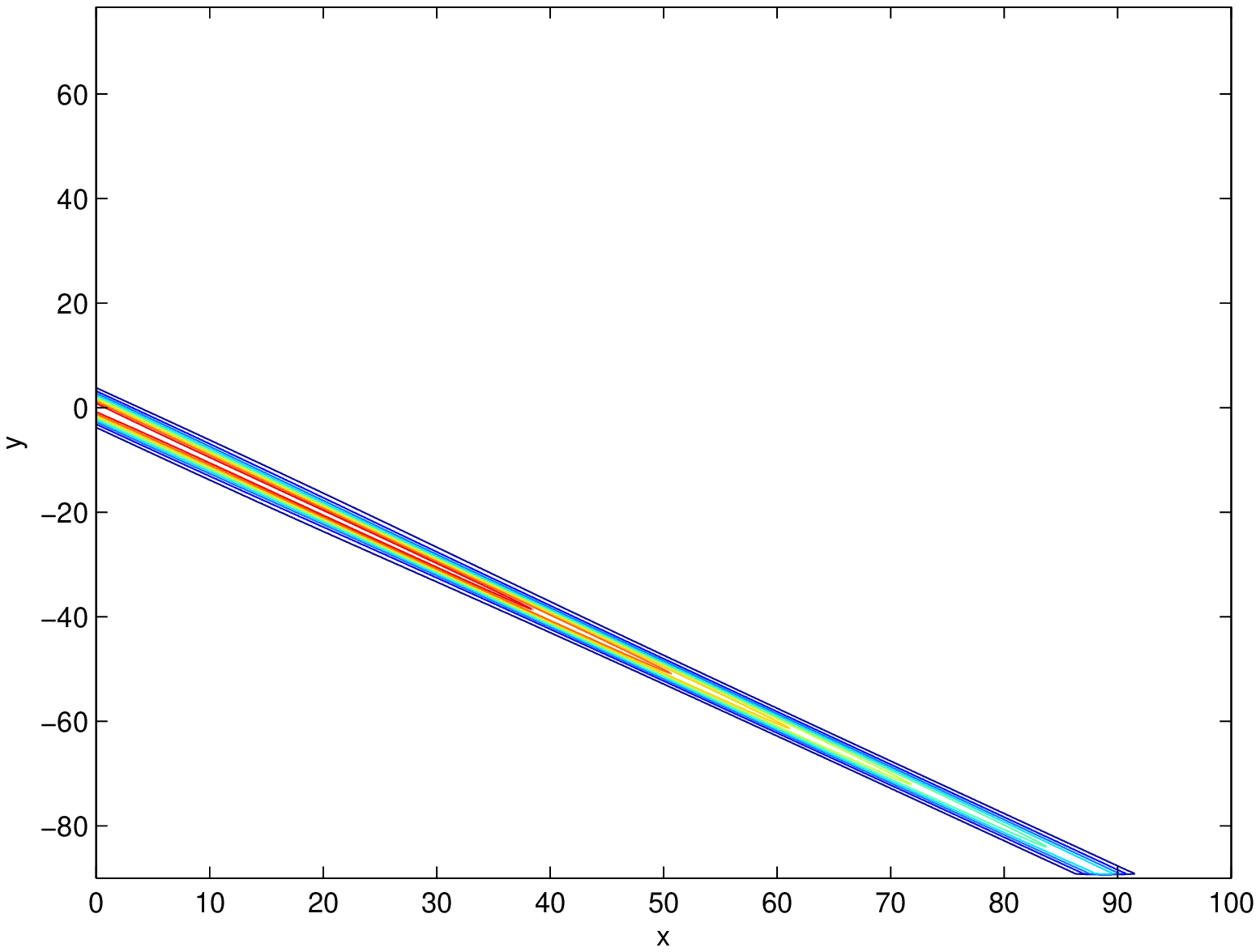, width=\linewidth, height=\linewidth}
\caption{\label{fig3-dvgce} First order scheme with $CFL=0.6,$ $\delta x=0.1,$ $\delta y =0.17.$ No focusing observed: the convergence of the scheme is poor.}
\end{minipage} \hfill
\begin{minipage}[b]{0.46\linewidth}
\centering\epsfig{figure= 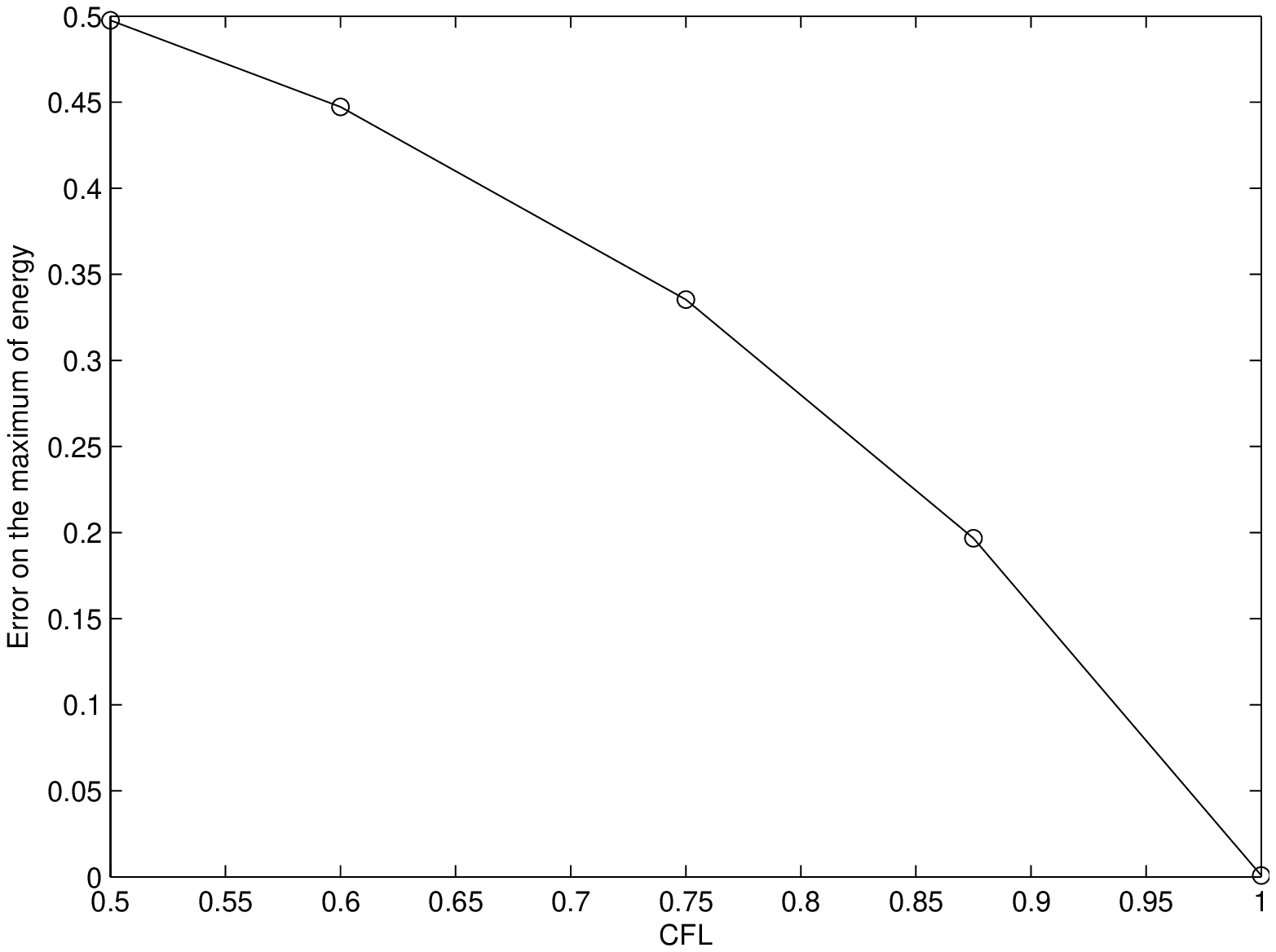, width=\linewidth, height=\linewidth}
\caption{\label{fig4-CFL_Max} First order scheme: error on the maximum of energy, as a function of $CFL$ (see Table \ref{tab:conv:ord1:CFL}).}
\end{minipage}
\end{figure}

If $CFL$ number decreases, the accuracy becomes bad  and even the focusing disappears: see Table  \ref{tab:conv:ord1:CFL}  and Figures  \ref{fig3-dvgce} and \ref{fig4-CFL_Max}. (Of course, if the CFL number is strictly larger than 1, the computed solution blows up).

\begin{table}
\begin{center}
\begin{tabular}{|l|ccccc|c|}
\hline
\bf{CFL} & $0.5$ & $0.6$  & $0.75$ & $0.875$ & 1 &  $\bf{1}$  \\ \hline
\bf{Error on energy} $\Sigma_{j,n}||u_j^n|^2-|u^{\rm ref,n}_j|^2|\delta x \delta y/|u^{\rm ref}|^{2} $ & 19 \% &	17\% &	14\% &	9\% &	2\%&	\bf{-}  \\ \hline
\bf{Focusing distance} & 43.1 & 49.1 &	55.6 &	48.0	& 59.9 &	\bf{59.7}  \\ \hline
\bf{Error on focusing distance} & 28\% &	18\% &	7\% &	19\%&	0.3\%	& \bf{-}   \\ \hline
\bf{Maximum of energy} & 1.08 &	1.18	& 1.42	& 1.72&	2.14 &	\bf{2.14}  \\ \hline
\bf{Error on the maximum of energy} 
 & 50\% &	45\%	&34\%	& 20\%	& 0.07\%&	\bf{-}  \\ \hline
\end{tabular}
\end{center}
\caption{Convergence of the first order scheme with cell size $\delta y=0.1$ and various $CFL.$  
The last column represents the fully converged reference case already seen $u^{ref}$ (with $\delta y=0.05$). We see that the focusing phenomenon is very poorly captured (huge error on the maximum of energy as soon as $CFL<1$).
\label{tab:conv:ord1:CFL}}
\end{table}

\

{\bf Convergence of the second order scheme}

\nopagebreak
We tested three different functions for the flux limiter: the first one is the Van Leer flux function defined by (\ref{eq:phi:VanLeer}), the second one is a convex combination of Lax-Wendroff and Beam-Warning flux limiter functions, defined by
\begin{equation}  \phi(\lambda)=\left\{
\begin{array}{ll}
0 & if\quad \lambda\leq 0 
\\
\lambda & if \quad 0\leq \lambda \leq 1
\\
1  & if \quad 1 \leq \lambda,
\end{array}
\right.
\end{equation}
the third one is the Superbee function defined by
\begin{equation}   \phi(\lambda)=\left\{
\begin{array}{ll}
0 & if\quad \lambda\leq 0 
\\
2\lambda & if \quad 0\leq \lambda \leq \frac{1}{2}
\\
1  & if \quad \frac{1}{2} \leq \lambda \leq 1
\\
\lambda  & if \quad 1\leq \lambda \leq 2 \\
2  & if \quad 2 \leq \lambda.
\end{array}
\right.
\end{equation}

\begin{figure}
\begin{minipage}[b]{0.46\linewidth}
\centering\epsfig{figure=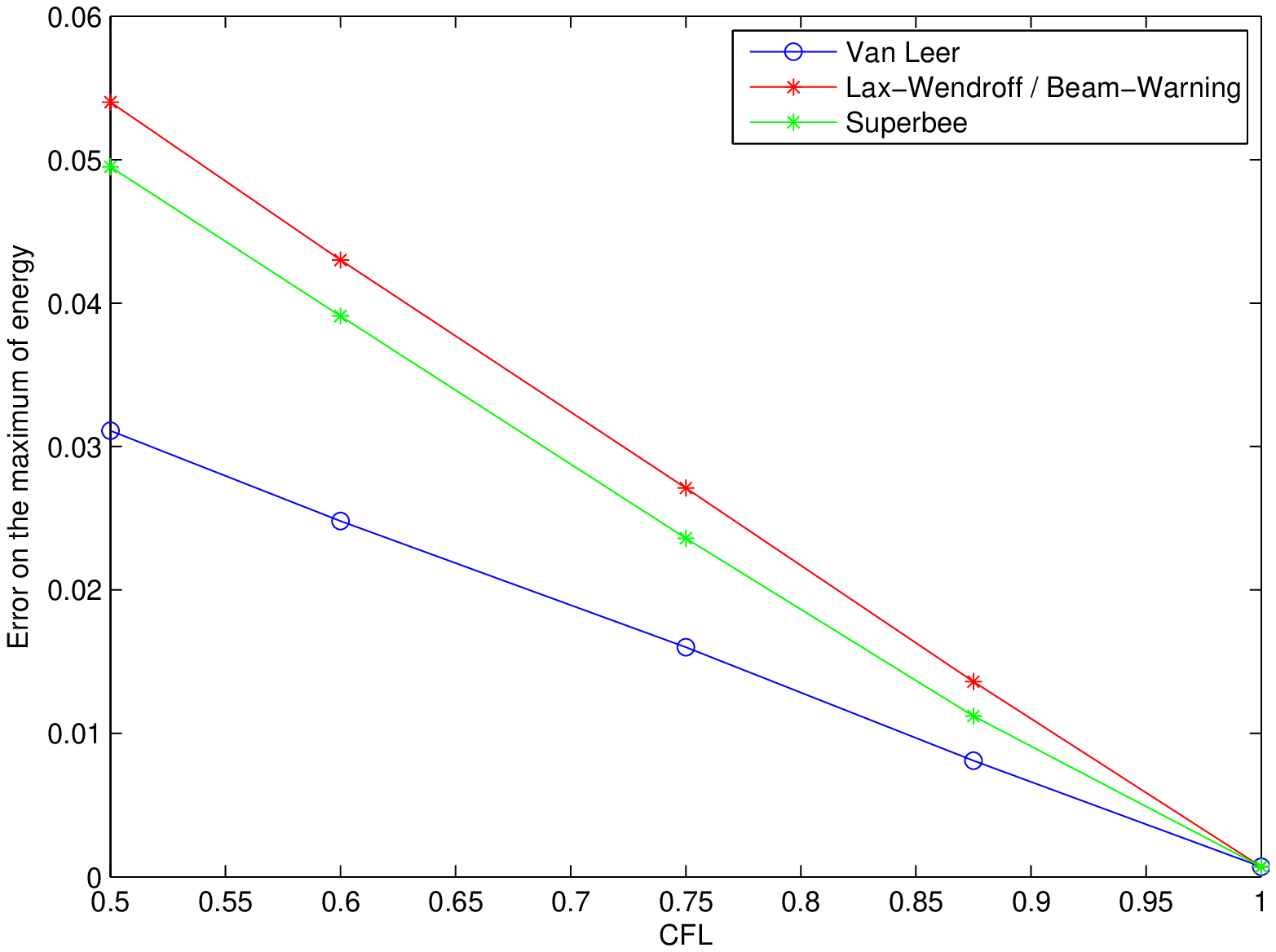, width=\linewidth, height=\linewidth }
\caption{\label{fig:Ord2Flux} Error on the maximum of energy as a function of $CFL,$ for $\delta y=0.1,$ for 3 different flux limiters.}
\end{minipage} \hfill
\begin{minipage}[b]{0.46\linewidth}
\centering\epsfig{figure=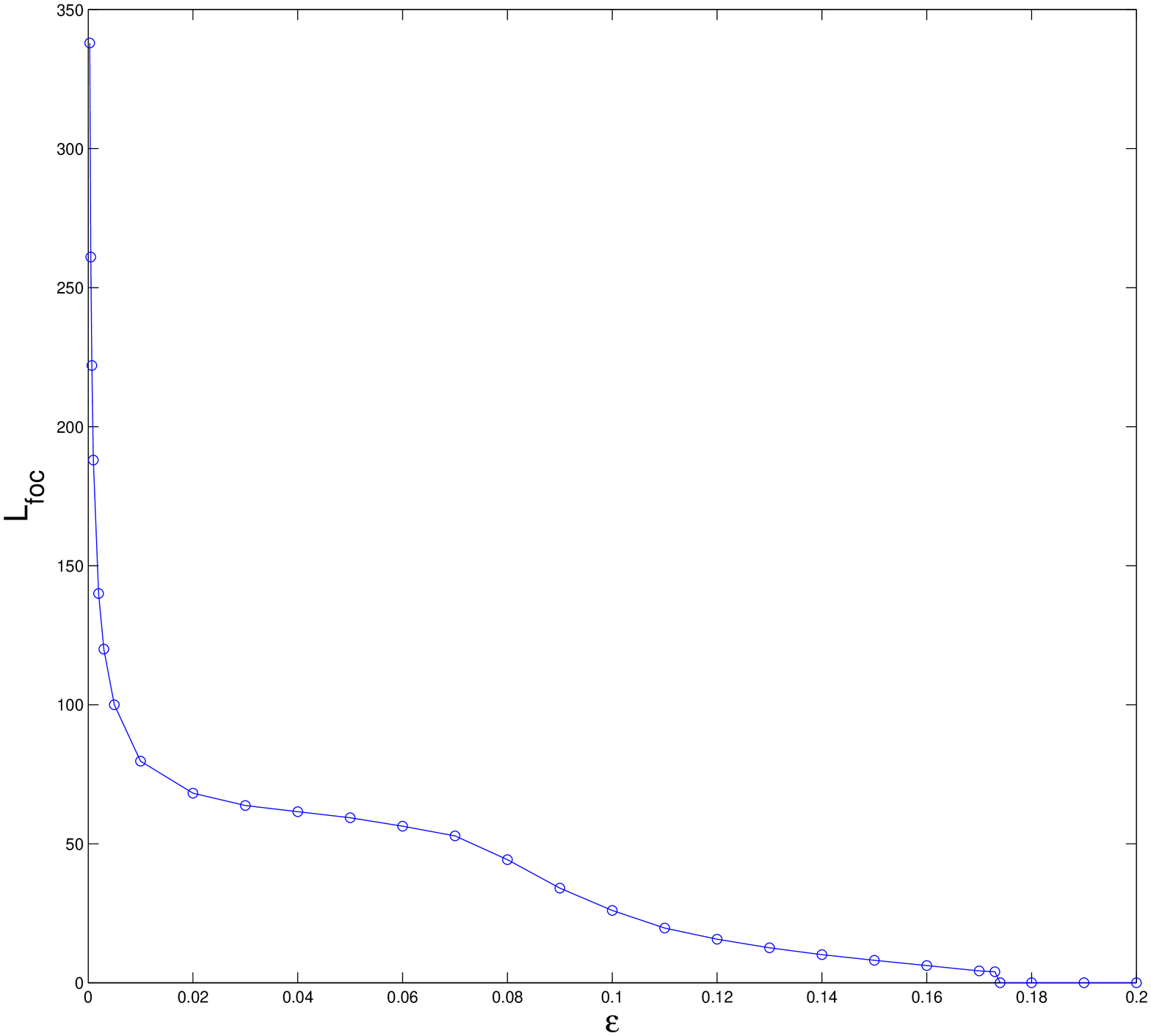,  width=\linewidth, height=\linewidth}
\caption{\label{fig:ep:Lfoc} Incidence of the variation of $\ep$ on the focusing distance (all other parameters as in the reference case, except $L_x$ and $L_y$).}
\end{minipage}
\end{figure}

We always apply these flux limiter functions at $\lambda=|u|^2$ and not at the real or imaginary part of the solution. As clearly shows Figure \ref{fig:Ord2Flux}, it appears that the Van Leer flux function is the one which gives the most accurate results. It is particularly clear in terms of the error on the maximum of energy : even for small CFL, its estimate is quite accurate contrarily to the first order scheme (for $CFL=0.5$ , the error is only about $3\%$ with second order scheme but about $50\%$ with first order one). 

\begin{figure}
\begin{minipage}[b]{0.46\linewidth}
\centering\epsfig{figure=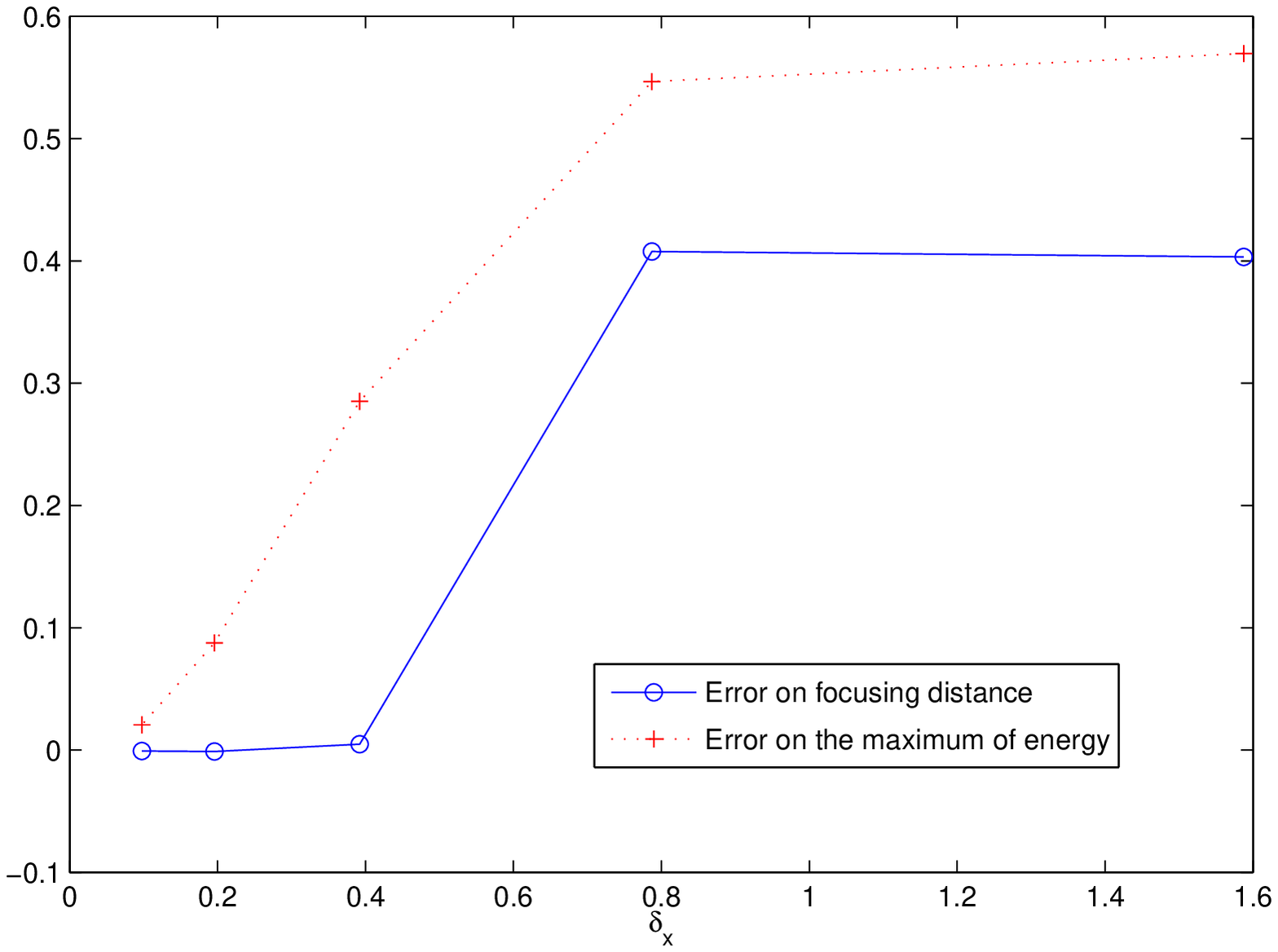, width=\linewidth, height=\linewidth }
\caption{\label{fig:ord2:CFL0.8} $CFL=0.8,$ second order scheme with Van Leer flux limiter: error on the focusing phenomenon as a function of the cell size $\delta x.$}
\end{minipage} \hfill
\begin{minipage}[b]{0.46\linewidth}
\centering\epsfig{figure=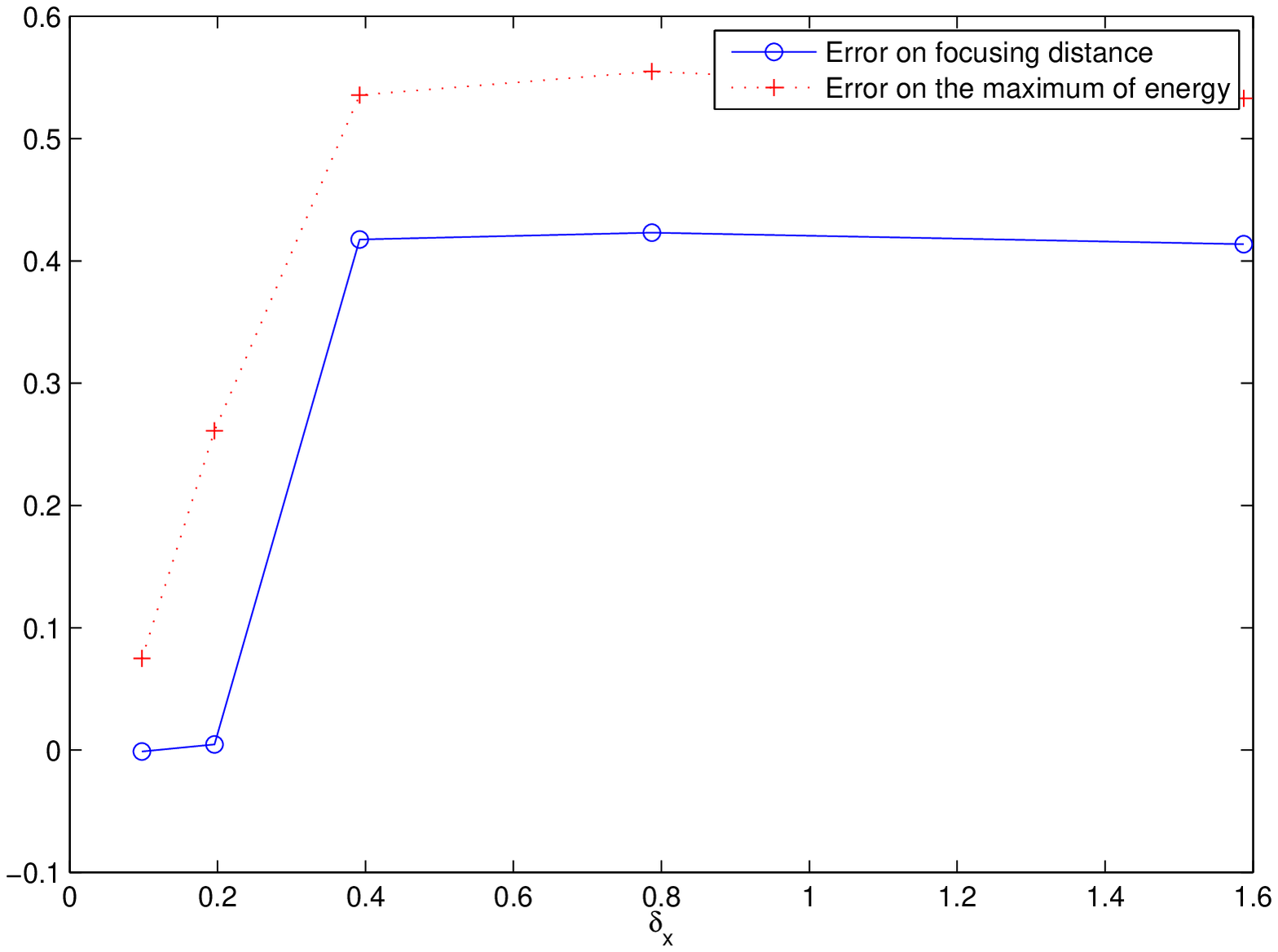,  width=\linewidth, height=\linewidth}
\caption{\label{fig:ord2:CLF0.6} $CFL=0.6,$ second order scheme with Van Leer flux limiter: error on the focusing phenomenon as a function of the cell size $\delta x.$}
\end{minipage}
\end{figure}

The smaller the CFL is, the more points are needed to get a correct approximation, as illustrates a comparison between Figures \ref{fig:ord2:CFL0.8} and \ref{fig:ord2:CLF0.6}. It is however performed even with $2^9$ points (that is, with $\delta x=0.2$) for $CFL=0.6$ for instance, contrarily to the scheme of order one, where no focusing at all is observed if $CFL=0.6$ even for $\delta x=0.1$ for instance (see Figure \ref{fig3-dvgce}).

\

{\bf Influence of the artificial boundary layer}

In the definition of the artificial absorbing layer $B$ given by (\ref{eq:defB}), we make $b$ and $\beta$ vary, with fixed cell sizes $\delta x=\delta y=0.2$ and all the other parameters given by the reference case. We look at the value of the total energy for each value of $b,\; \beta$  (the reference values being  $b=0.1,$ $\beta=50.$) The results are given in Table \ref{tab:B}. We check that the sensitivity to the exact values of these coeficients is very weak; but it is crucial to have $b \neq 0 ,$  elseif spurious reflexions may appear on the boundaries.

\begin{table}
\begin{center}
\begin{tabular}{|l|cccc|}
\hline
& $\bf{\beta=10}$ &$\bf{\beta=30}$ &$\bf{\beta=50}$ &$\bf{\beta=100}$ \\ \hline
{\bf b=0} &   29\%  &  29\%  &  29\% &    29\% \\ \hline
{\bf b=0.1} & 0.08\%  &  0.02\%  &  {\bf 0} &    0.02\% \\ \hline
{\bf b=0.2} &  0.03\% &  0.03\% &  0.05\% &  0.07\% \\ \hline
{\bf b=0.5} &  0.08\% &  0.14\% &   0.15\% &   0.16\% \\ \hline
{\bf b=1} &    0.19\% &   0.22\% &   0.23\% &   0.23\% \\ \hline
\end{tabular}
\end{center}
\caption{Incidence of the variation of the boundary layer $B$ on the difference between the total energy of each case and the one of the reference case ($b=0.1$ and $\eta=50$): $\Sigma_{j,n}||u_j^n|^2-|u^{\rm ref,n}_j|^2|\delta x \delta y/|u^{\rm ref}|^{2}.$  The results of this table show that the influence is negligible, as soon as $b$ is not zero. \label{tab:B}}
\end{table}

\subsubsection{ Variation of several parameters}

\

$\bullet$ {\bf Variation of the absorption coefficient}
\begin{table}
\begin{center}
\begin{tabular}{|l|ccccccc|}
\hline
& $\bf{\frac{\nu_0}{\nu}=0}$ & $\bf{\frac{\nu_0}{\nu}=0.1}$ & $\bf{\frac{\nu_0}{\nu}=0.3}$  &$\bf{\frac{\nu_0}{\nu}=0.5}$ &$\bf{\frac{\nu_0}{\nu}=0.7}$ & $\bf{\frac{\nu_0}{\nu}=0.9}$ & $\bf{\frac{\nu_0}{\nu}=1}$ \\ & & & & & & & \\ \hline
{\bf reference case:} & & & & & & & \\ {\bf $\nu=10^{-3},$ $\alpha=0.05$} & 0.3\% &   0.2\% &    0.1\% &         - &   0.1\% &    0.2\% &    0.3\%     \\ \hline
{\bf $\nu=10^{-3},$ $\alpha=0.5$} & 6.2\% &   5.0\% &     2.5\%   &-&    2.5\% &    5.0\% &   6.2\% 
\\ \hline
{\bf $\nu=10^{-2},$ $\alpha=0.05$} &  0.5\% &    0.4\% &    0.2\% &         - &   0.2\%   & 0.4\% &    0.5\%  
\\ \hline
{\bf $\nu=10^{-2},$ $\alpha=0.5$} & 8.9\% &    7.2\% &    3.6\% &         - &    3.7\% &    7.4\% &    9.3\% 
\\ \hline

\end{tabular}
\end{center}
\caption{Influence of the repartition between $\nu_0$ and $\nu_1$ in different cases:  percentage of error on total energy, defined by $\Sigma_{j,n}||u_j^n|^2-|u^{\rm ref,n}_j|^2|\delta x \delta y/|u^{\rm ref}|^{2}.$ \label{tab:nu}}
\end{table}

The numerical scheme can also be used with no absorption ($\nu=0$), it still works and give good results. The repartition of $\nu_0$ and $\nu_1$ changes very little the solution, as shows Table \ref{tab:nu}. In each case, the reference is taken for $\nu_0=\nu_1=\frac{\nu}{2}.$ The table shows the results only for the comparison on the total energy; indeed, the focusing distance remains completely unchanged in any case, and the maximum of energy changes by less than $0.3\%$ in the worst case.

When the absorption coefficient is larger, the problem is easier to solve since the laser energy decreases when  $x$ increases: for instance in the reference case, if we set $\nu=10^{-2}$ instead of $\nu=10^{-3},$ the ray is rapidly totally absorbed, and no focusing is observed. 

The influence of the repartition between $\nu_0$ and $\nu_1$ increases with $\alpha,$ as shows Table \ref{tab:nu}.               

\

$\bullet$ {\bf Variation of the incidence angle}
   
To test whether the scheme is accurate for various angles, we make it vary from $5^0$ to $70^0,$ all the other parameters being constant: see Table \ref{tab:angle}.
 We check  that the indicators for the focusing distance and the maximum of energy are  well  estimated, since they depend very few on the incidence angle.

\begin{table}
\begin{center}
\begin{tabular}{|l|ccccc|}
\hline
{\bf Incidence angle}              & $\bf{5^0}$ &$\bf{30^o}$ &$\bf{45^0}$&$\bf{60^0}$ & $\bf{70^0}$ \\ \hline
${\bf \delta x }$                  &  0.23      &  0.16      & 0.2       &    0.16   & 0.02  \\ \hline
${\bf \delta y}$                   & 0.02       &  0.1       & 0.2       &    0.27   & 0.06  \\ \hline
{\bf Maximum of energy}                &  2.17      &  2.16     &  2.13     &  2.10     & 1.99  \\ \hline
{\bf Error on the maximum of energy}   & 1.5\%      &  0.8\%    &   0.43\%  &   2.2\%   & 9.7\% \\ \hline
{\bf Focusing distance}             &    59.2    &   59.7    &   59.35   &   59.9    &60.2  \\ \hline
{\bf Error on the Focusing distance}&    0.9\%   &   0.01\%   &   0.6\%  &   0.34\%  & 0.96\% \\ \hline

\end{tabular}
\end{center}
\caption{Variation of the incidence angle: influence on the focusing distance and on the maximum of energy. As usual, the errors refer to the fully-converged reference case.
\label{tab:angle}}
\end{table}

\

$\bullet$ {\bf Variation of $\ep$}

 If all other coefficients are fixed, the larger $\ep $ becomes,  the more important the diffusion phenomenon is (and the larger the domain must be  to obtain a converging solution), and, in the nonlinear case, the smaller the focusing distance becomes. A limit value of $\ep$ is experienced, above which no focusing phenomenon (for the nonlinear equation) is observed. In our reference case for instance, the limit is around $\ep=0.17,$ see Figure \ref{fig:ep:Lfoc}, but this limit depends of course on all parameters, especially $\alpha$ and $\nu.$

From a physical point of view, all our asymptotic analysis is built on the assumption $\ep=o(1):$ else, our equation is no more a valid approximation of the envelope of Helmholtz equation, given by (\ref{base}). Hence, we have to assume $\ep << 1:$ larger values are meaningless.  


\

$\bullet$ {\bf Variation of $\alpha.$}

The parameter $\alpha$ represents a nonlinear effect, and induces autofocusing and filamentation of the beam.
The larger it is, the more accurate the focusing phenomenon becomes, as illustrated in Figure \ref{fig:alpha:Max}. 

It could be interesting to evaluate the value of $\alpha$ for which a focusing  phenomenon appears:
in our reference case, it is for $\alpha \ge 0.02$. 
On the other hand, one may check that if $\alpha$ is large
 enough, several focusing points appear and a breaking of the beam occurs (see Figure \ref{fig:alpha:increase:ray}). This phenomenon depends of course also on the absorption coefficient $\nu$ and on the diffusion coefficient $\ep.$  

\subsubsection{Remark on artificial damping}
\begin{figure}
\begin{minipage}[b]{0.46\linewidth}
\centering\epsfig{figure=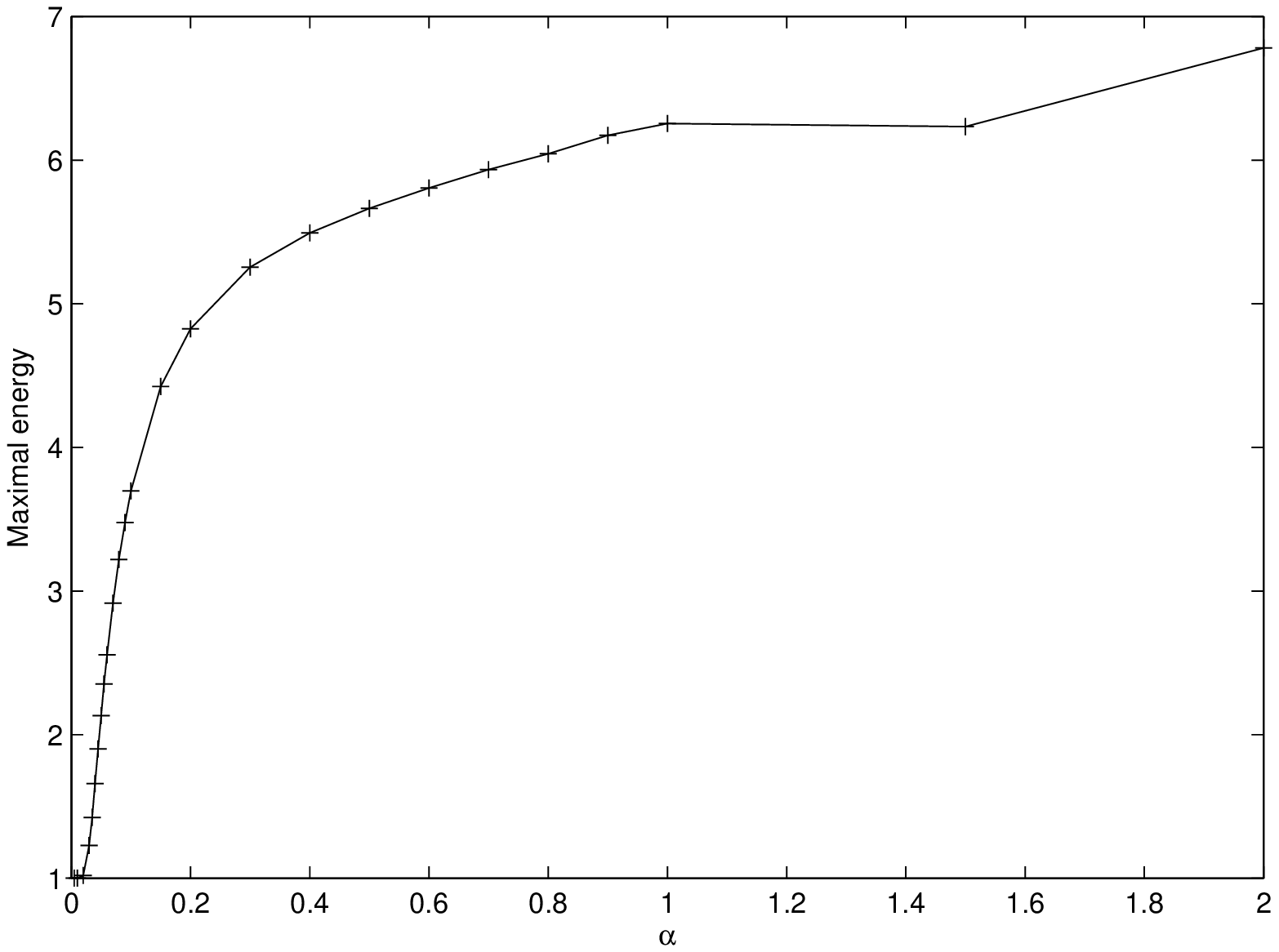,width=\linewidth,height=\linewidth}
\caption{\label{fig:alpha:Max} Influence of $\alpha$ on the maximum of energy (obtained in the focusing phenomenon). Standard hypothesis.The autofocusing, which is a nonlinear effect, is more significant when $\alpha$ increases.} 
\end{minipage} \hfill
\begin{minipage}[b]{0.46\linewidth}
\centering\epsfig{figure=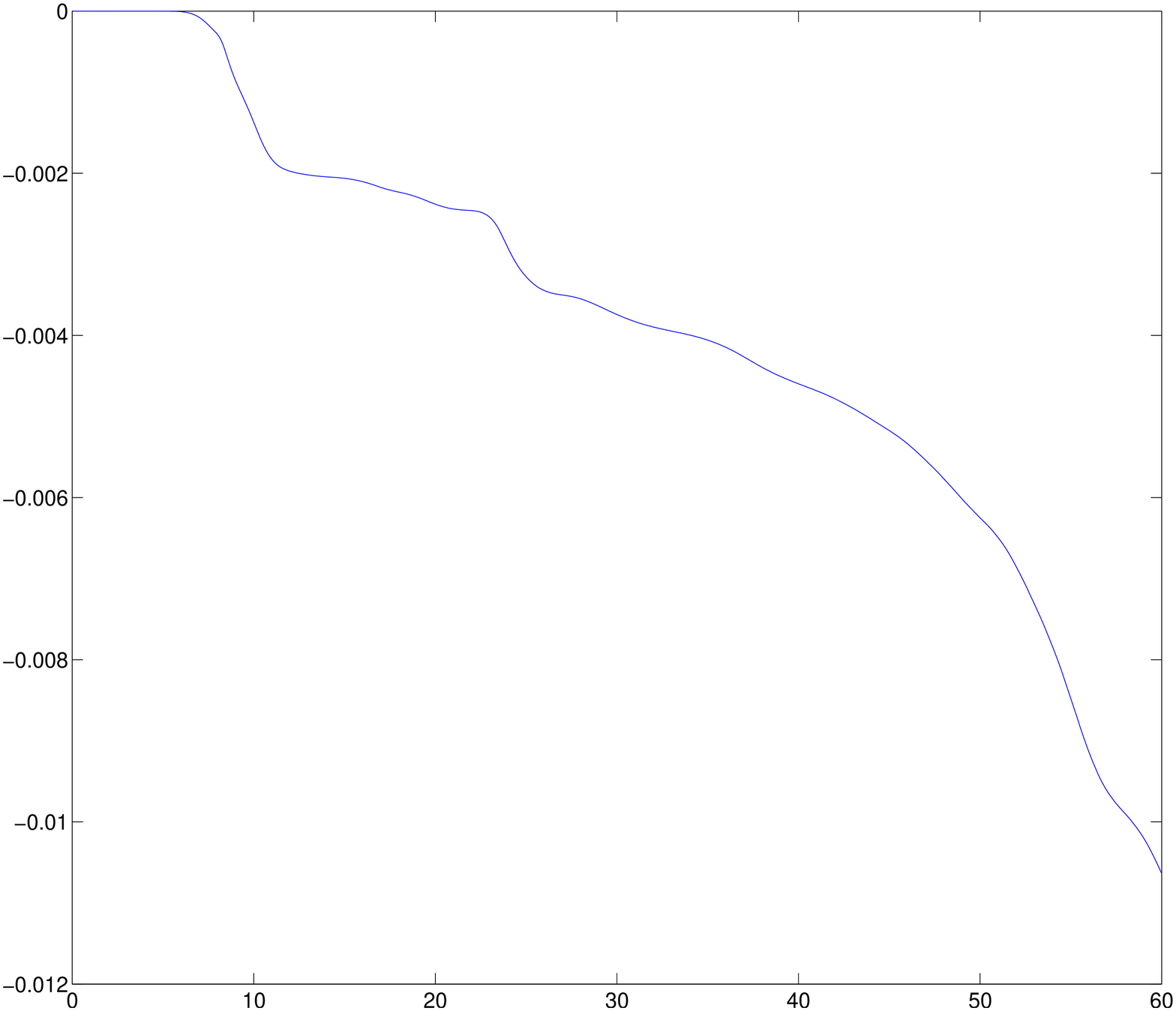, width=\linewidth, height=\linewidth}
\caption{$\alpha=1.5,$ $\nu=0:$ we define the energy $E^n=\Sigma_j |u^n_j|^2 .\delta y.$ This picture shows $(E^n-E^0)/E^0$ as a function of $x^n=n\delta x:$ the energy $E^n$ decreases by less than $2\%$ during the whole trajectory.\label{fig:nonu:energy}}
\end{minipage} \hfill
\end{figure}

We wish to check now that there is no artificial damping due to the numerical scheme; in other words, that in the second stage the decrease of the $l^2-$ norm of the solution has the right value. Using the notations of Section \ref{subsec:scheme:prop}, this right value is given by the equality:
$$||u^{n+1}||_{l^2}=e^{-2\nu_1\f{\delta x}{kx}}||u^{n\#}||_{l^2}.$$
 Going back to Equation (\ref{eq:step2}), we can write it under the form (assuming no artificial boundary layer: $B_j=0$)
$$u^{n+1}_j=\frac{1-a-ib}{1+a+ib}u^{n\#}_{\theta_j},$$ 
where we set $a= \frac{\delta x}{2k_x} \nu_{1,j}^n$ and $b=\frac{\delta x}{2k_x}\mu_j.$ 
Since the characteristic value of the coefficient $a$ is $10^{-4}$ (or smaller) and, in the worst case, the characteristic value of $\mu$ is in the order of $1,$ so that we can choose $\frac{\delta x}{2k_x}$ to have $b$ small, we see that
$$|\frac{1-a-ib}{1+a+ib}|^2=1-4a\frac{1}{1+b^2}+o(a^2),$$
which is very close to the right value $e^{-4a}=1-4a+o(a^2).$ The only damping may then come from the fact that $\sum\limits_j |u^{n\#}_{\theta_j}|^2$ may be significantly smaller than $\sum\limits_j |u^{n\#}_j|^2,$ due to a large difference between $u^{n\#}_j$ and $u^{n\#}_{j-1}.$ 
To check this numerically, we test the case $\nu=0:$ Figure \ref{fig:nonu:energy} shows that even in a difficult case with a large $\alpha=1.5,$ the global energy $||u^n||^2_{l^2}$ is conserved.

\subsubsection {Two-ray model}

\begin{figure}
\begin{minipage}[b]{0.46\linewidth}
\centering\epsfig{figure=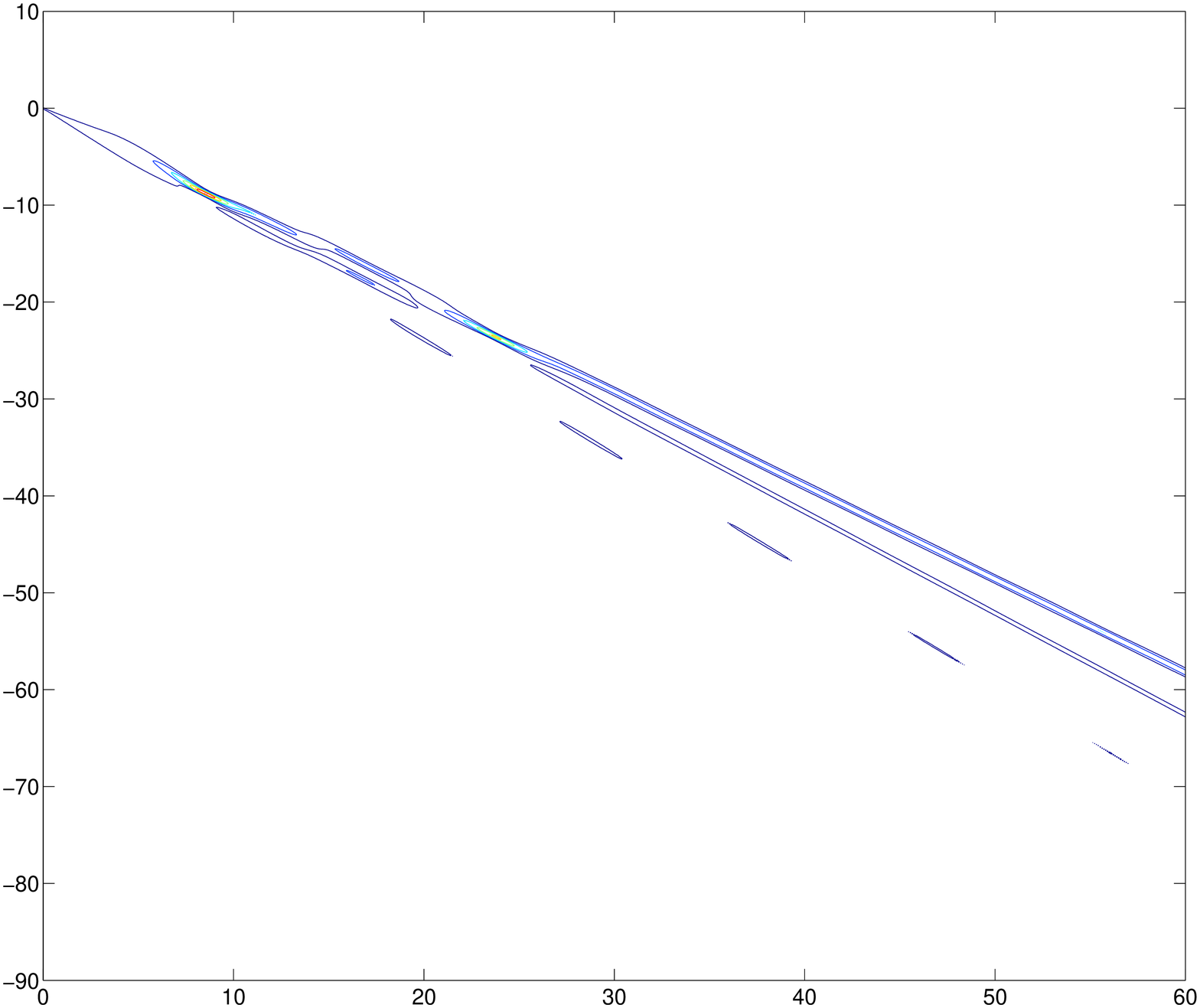,width=\linewidth,height=\linewidth}
\caption{\label{fig:alpha:increase:ray} $\alpha=1.5,$ $\nu=0:$ high focusing. One observes a breaking of the beam in three sub-beams.}
 \end{minipage} \hfill
\begin{minipage}[b]{0.46\linewidth}
\centering\epsfig{figure=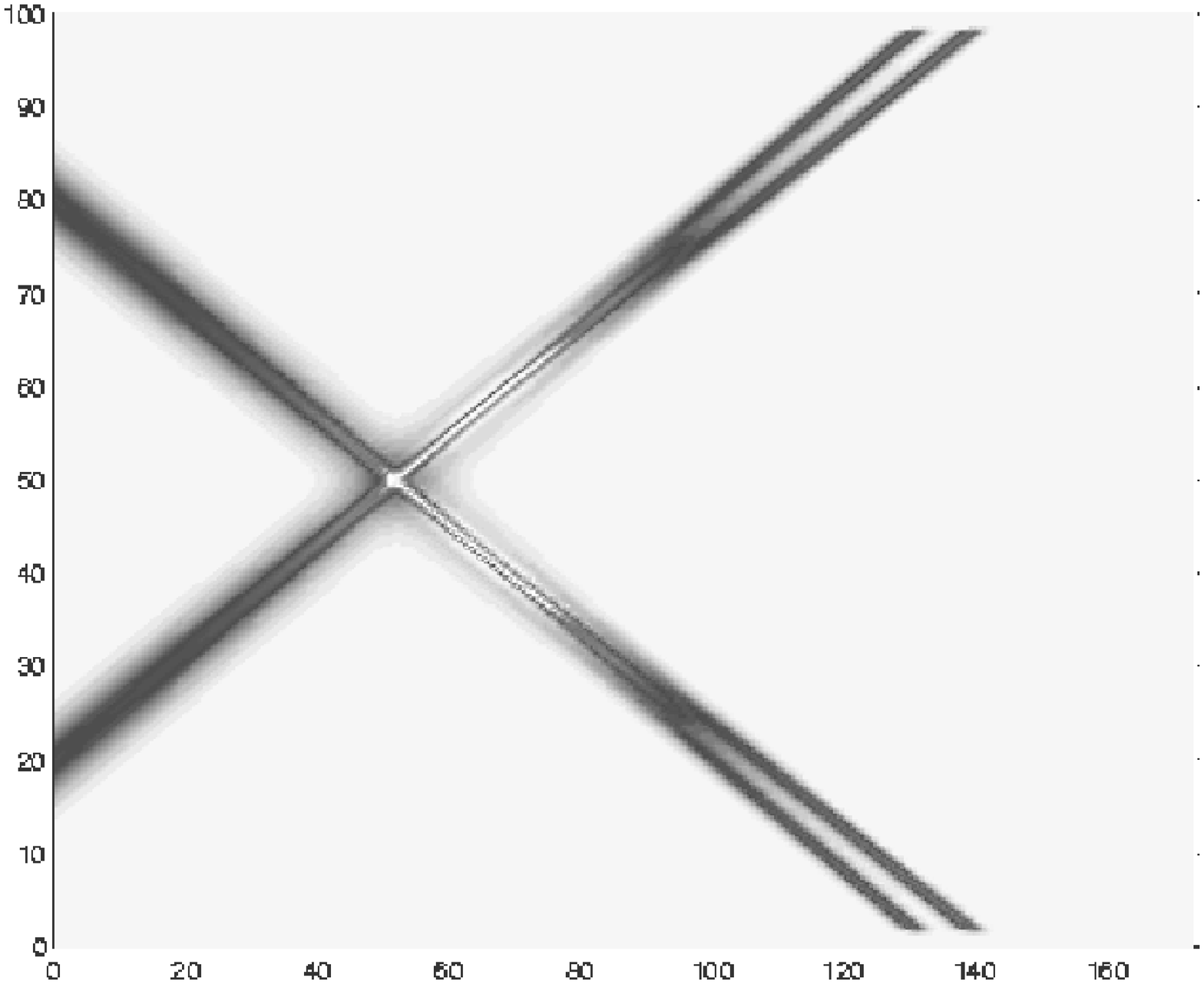, width=\linewidth, height=\linewidth}
\caption{  2 beams crossing with incidence angles $\pm 30^0,$ $\alpha=0.05,$ and $L=5$ for the initial gaussian functions.}
\label{fig:2rays}
\end{minipage} \hfill
\end{figure}

We have also performed computations for the two-ray model which is described above at Section \ref{subsubsec:2ray} using two functions $u^1$ and $ u^2$;
an illustration is given by Figure \ref{fig:2rays}. The interaction between the rays is only given by the nonlinear term $f(w)$ with $w^2=|u^1|^2+|u^2|^2$ as above. To analyse its exact influence, one can compare the result given by the previous model with the two-ray interaction  and the result given by a simple  superposition of two independant rays (obtained with the one-ray model). One may see then that the energy becomes larger with the two-ray interaction: on the case of Figure \ref{fig:2rays} for instance, $Max(|u^1|^2+|u^2|^2)=12.3$ instead of $10.6$ if the rays do not interact.

\section {Extension to a Time-Dependent Interaction Model}\label{sec:tilted}

We now address  a model where a tilted paraxial 
 equation is coupled with a hydrodynamic model in order to
study  filamentation.
 Under the hypothesis of a small incidence angle, this model has been extensively used by physicists for a long time and it is also addressed in \cite{para-3},\cite{fcs},\cite{dorr} for example and the references therein (for a derivation of this model, see \cite{sen} for example).

\subsection{ The Model and the Numerical Method}

\textbf{Modeling of the plasma.}

By taking  the critical
density  (depending only on the laser wave length) as a reference density, one defines a non-dimension electron
density  $N=N(t,\mathbf{x})$ ; so the plasma may be characterized only by  this quantity, the plasma velocity $\mathbf{U}=\mathbf{U}(t,\mathbf{x})$ and the electron density $T_e(t,\mathbf{x})$.

   Then, the simplest model is the 
following one.  The  pressure $P=P(N, T_e)$
is assumed to be a smooth function of the density $N$ and of the electron temperature
$T_e$ (which is assumed to be a very smooth fixed function of the position $\mathbf{x}$ ), for example $P(N, T_e)$ may be  the sum of two terms equal to $N^3 $ and $ N T_e$ up to multiplicative constants. Then one considers the following
barotropic Euler system:
\begin{eqnarray}
\frac{\partial }{\partial t}N+\nabla (N\mathbf{U}) &=&0,  \label{moli} \\
\frac{\partial }{\partial t}(N\mathbf{U})+\nabla (N\mathbf{UU})+\nabla
(P(N, T_e)) &=&-N\gamma _{p}\nabla |\Psi |^{2}.  \label{molii}
\end{eqnarray}

The term $\gamma _{p}\nabla |\Psi|^{2}$ corresponds to a ponderomotive force
due to a laser pressure (the coefficient $\gamma _{p}$ is a constant
depending only on the ion species).

\textbf{Modeling of the laser beam. } 

The laser
field $\Psi =\Psi (t,\mathbf{x})$ is a solution to the following frequency wave equation 
(which is of Schr\"{o}dinger type):
\begin{equation}
2i\frac{1}{c}\frac{\partial }{\partial t}\Psi+\frac{1}{k_{0}}\Delta
\Psi+k_{0}(1-N) \Psi+i\nu^{\diamond} \Psi=0,  
\label{wave}
\end{equation}
where the real coefficient $\nu^{\diamond} $ is related
 to the absorption of the laser intensity by the plasma
and  $c$  the light speed.

 Assume that the mean value of
the plasma density is quite constant and denoted  $N_{m}$, so we set:
$$ 
N(\f{x})=N_{m} +\delta N(\f{x}),
$$
where $\delta N$ is small with respect to $1.$ Then one can make the paraxial approximation ; that is to say the laser beam is now characterized by the space and time
 envelope of the electric field $\mathcal{U}=\mathcal{U}(t,\mathbf{x})$ and we set:
$$
 \Psi(t,\f{x})=\mathcal{U}(t,\f{x})
 e^{ik_{0}\f{K}.\f{x} }, \qquad {\rm where } \; \f{K}=
 \sqrt{1-N_{m} }\f{k}.
$$ 
 Therefore, if one sets $\epsilon= \frac{1}{k_{0} \sqrt{1-N_{m}}} ,$ by the same procedure as mentioned in the introduction,  one checks that $\mathcal{U}$  
 satisfies:
\begin{equation}
\sqrt{1-N_{m}} (i\mathbf{k.}\nabla \mathcal{U} +\frac{\epsilon}{2} \Delta _{\bot}^{k} \mathcal{U} )
 +i\frac{\nu^{\diamond}} {2}    \mathcal{U}  -\frac{k_{0} \delta N} {2 }  \mathcal{U} + i\frac1{c }\frac{\partial \mathcal{U} }{\partial t}=0.  \label{parax}
\end{equation}

 It is necessary to supplement  equation
(\ref{parax}) with  the same boundary condition 
as in the model of section 1  (and with an initial condition).
  
{\bf Numerical method.}

We consider a mesh of finite difference type as above. 
The numerical treatment of the barotropic Euler system (\ref{moli})(\ref{molii}) is a classical one, we have chosen a \emph{Lagrange-Euler} method, see \cite{fcs} for details.
 To deal with (\ref{parax}), according to the  large value of the speed of light, one must perform a time inplicit discretization. So at each time step, one solves firstly the  Euler system with a ponderomotive force evaluated with the previous value of $|\mathcal{U}|^2. $ Secondly, using the obtained values of $N$ and of $\delta N ,$ one has to solve (\ref{parax}) ;  if $u^{ini}$ and $u$ denote the values of the field $\mathcal{U}$ at the beginning and the end of time step, one searches $u$  solution to:
\begin{equation}
i \mathbf{k.}\nabla u
+i \nu u +\frac{\epsilon}{2} (\Delta _{\bot
}^{k}u)-\mu  u=
\frac{i}{c \sqrt{1-N_{m}} }\frac{u^{ini}}{\delta t},
\label{parr}
\end{equation}
 where we have set:
$$ \mu= \frac{k_{0}\delta N} {2 \sqrt{1-N_{m}} } , \qquad \nu=\frac{1}{c \sqrt{1-N_{m}}  }\frac{1}{\delta t}+\frac1{2\sqrt{1-N_{m}} } \nu^{\diamond}.
$$
That is exactly the equation studied in section 3, but a right hand side term has been added.
So the numerical method is the same as described above ; the only modification is the adding of the right hand side term  in the transport stage. 
Notice that the index of refraction $(1-N)$ is equal to $(1-2\epsilon \mu )(1-N_m).$

>From a practical point of view, the numerical method for (\ref{parax}) has been implemented in a parallel way in the {\it{HERA}} plateform for plasma hydrodynamics in 2D and in 3D; the parallel solver and the domain decomposition techniques are the same as the ones detailed in \cite{fcs}.

\subsection {Numerical Results}

\begin{figure}[t]
\begin{minipage}[b]{0.46\linewidth}
\includegraphics[height=4.2cm]{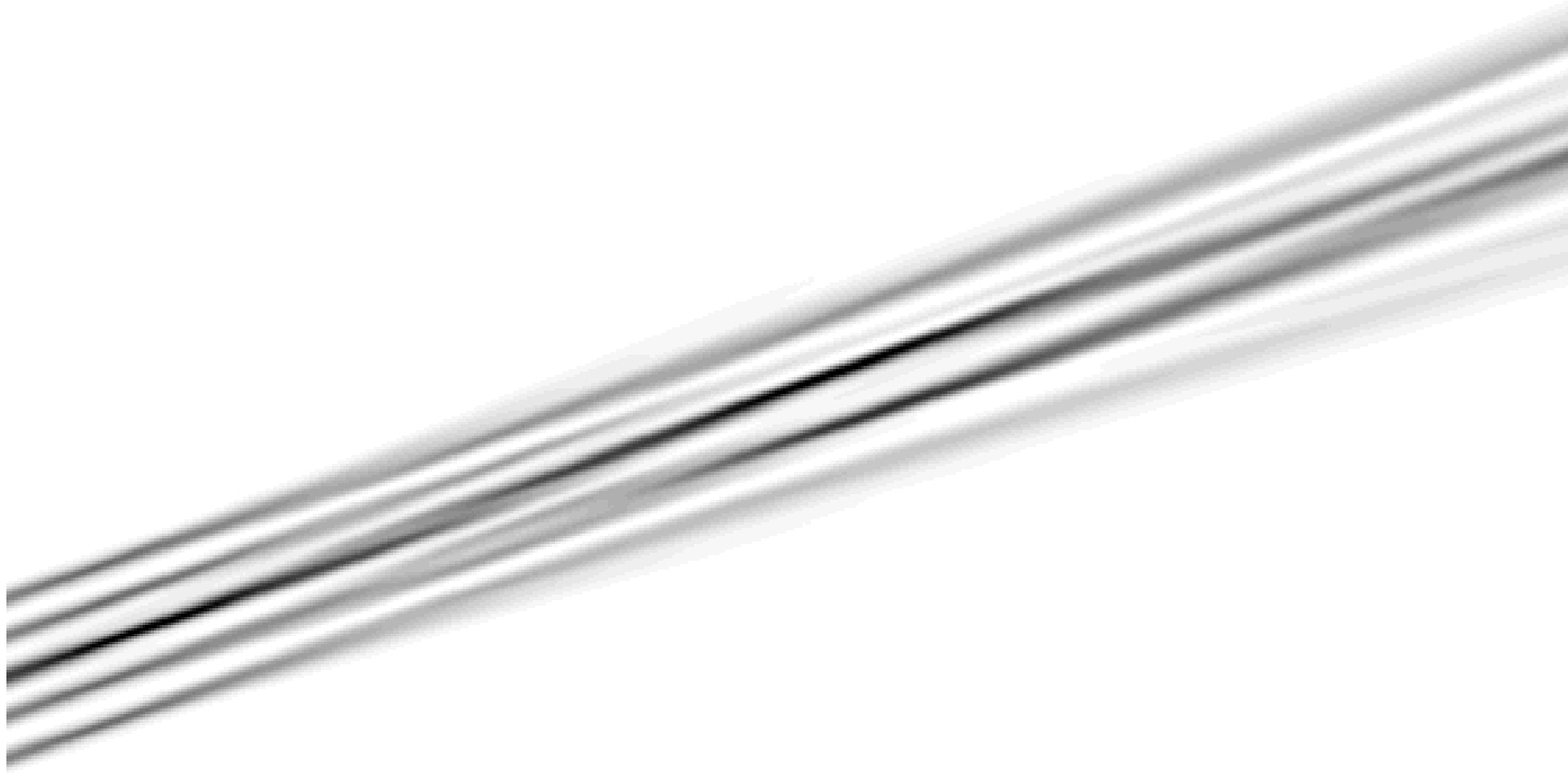} 
\end{minipage} \hfill
\begin{minipage}[b]{0.46\linewidth}
\includegraphics[height=4.2cm]{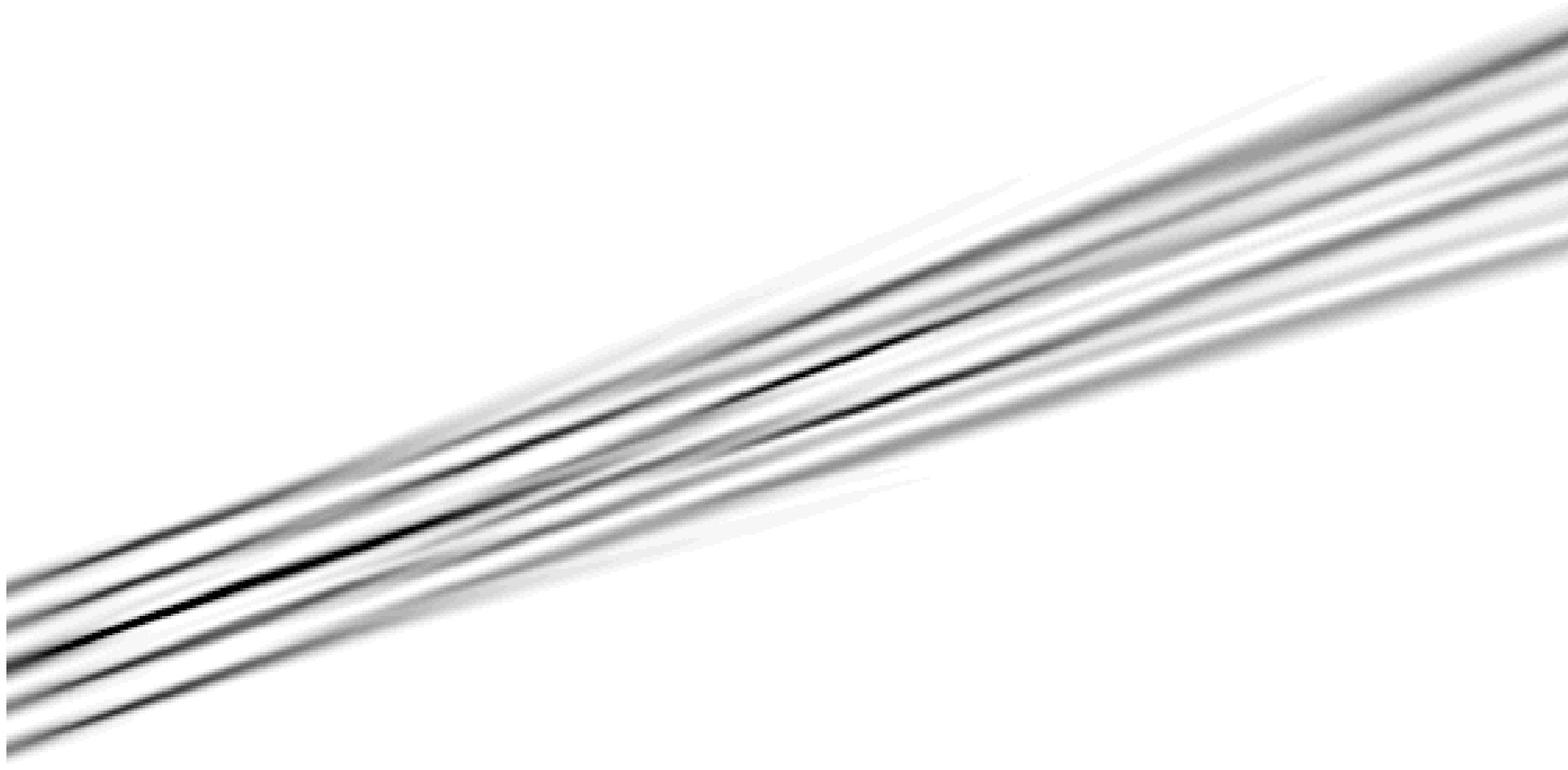} 
\end{minipage}
\begin{minipage}[b]{0.46\linewidth}
\includegraphics[height=4.2cm]{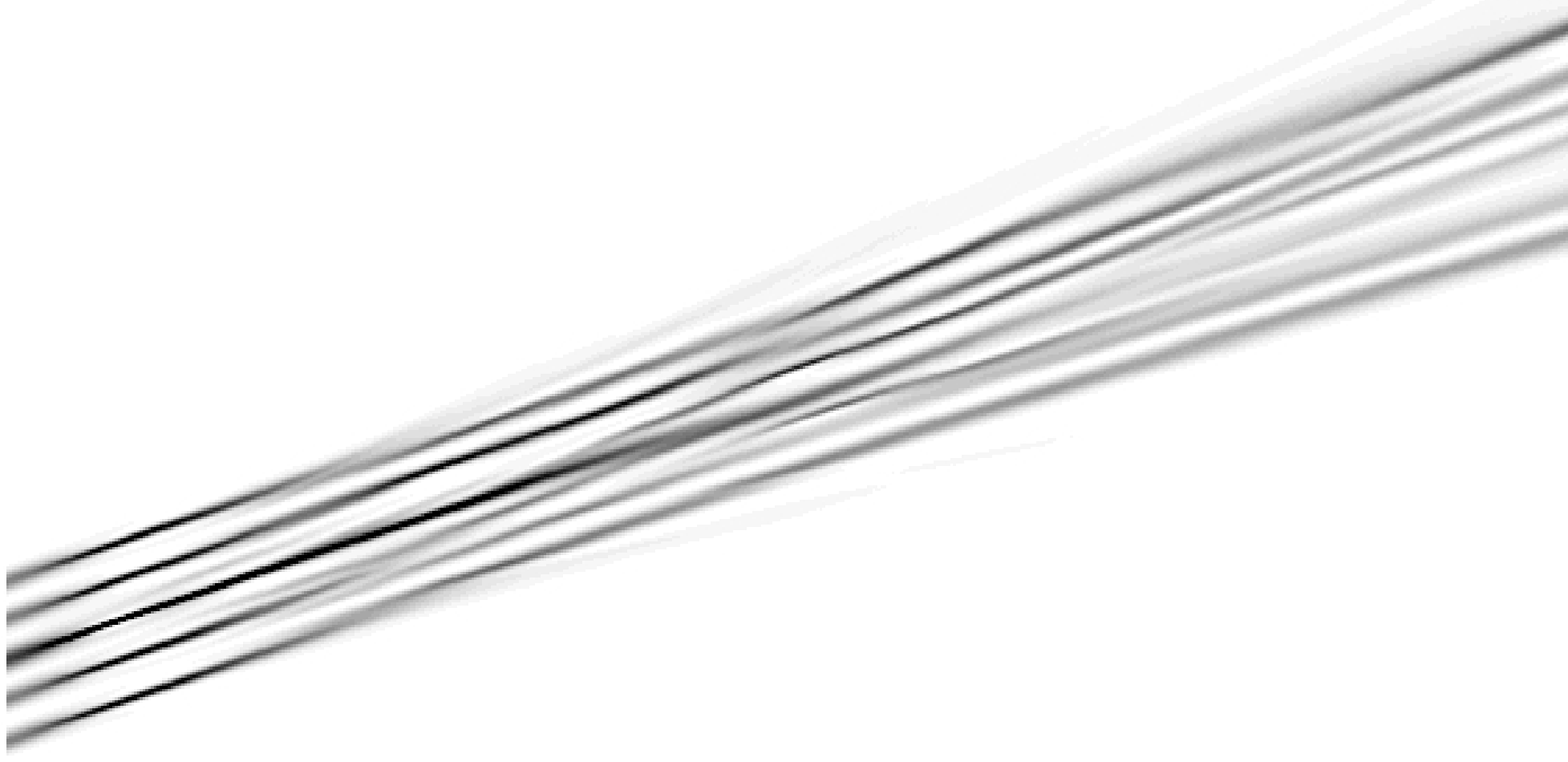} 
\end{minipage} \hfill
\begin{minipage}[b]{0.46\linewidth}
\includegraphics[height=4.2cm]{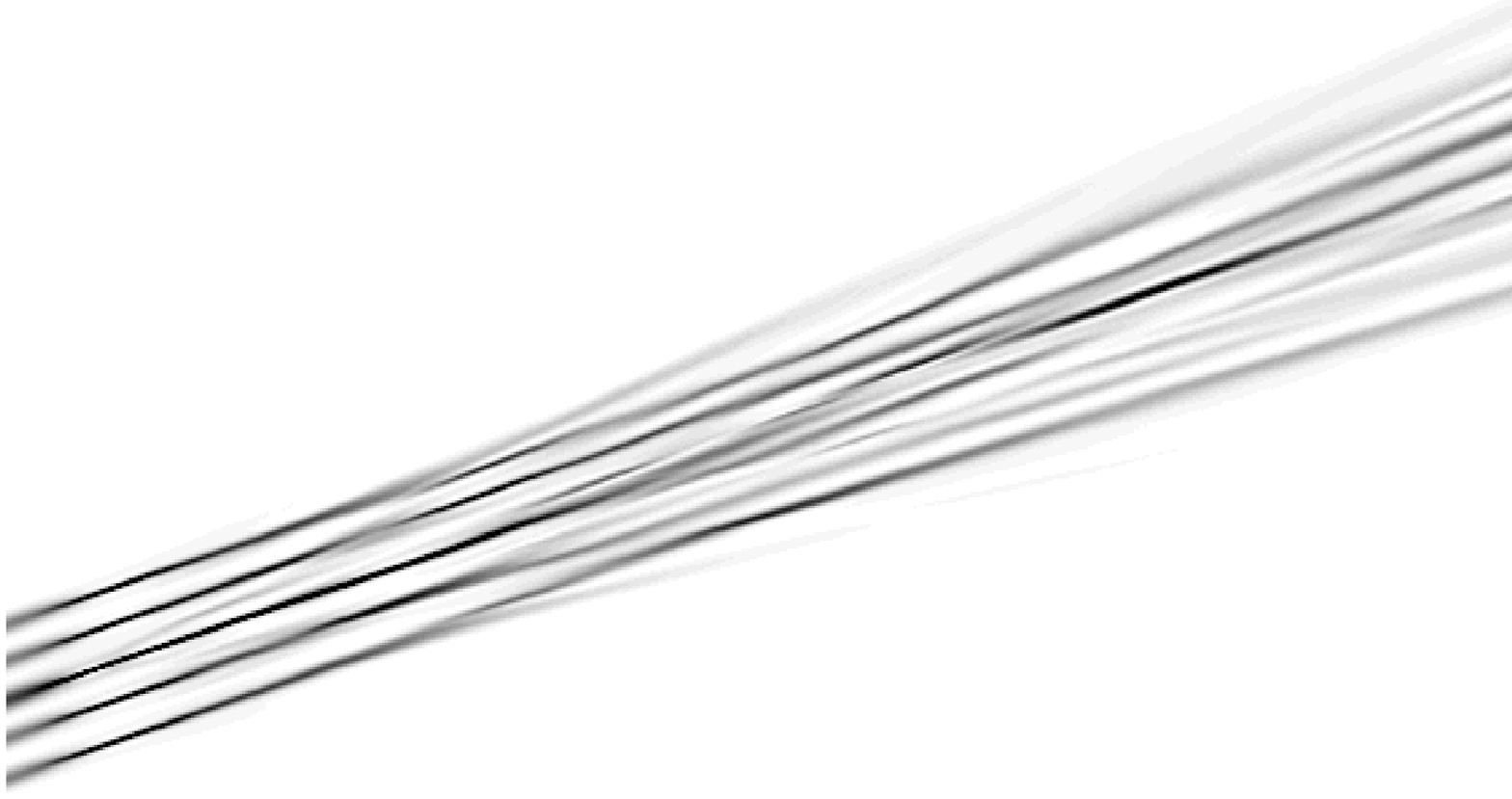} 
\end{minipage}
\caption{{\small Snapshot of the laser intensity at the time 2.6 ps, 3.9 ps, 5.3 ps and 6.6 ps ( from the top-left to the bottom-right).\label{beau}}}
\end{figure}

Recall that from a practical point of view, in the 
transverse profile of a laser beam, one distinguishes a lot of small hot spots, called {\it speckles}, 
whose intensity is very large compared to the mean intensity of the beam.
The shape of each individual speckle is a Gaussian function whose width is about
a few micrometers. We present here the results of a 2D numerical  simulation. One addresses a simulation box which is 600 $\mu$m  long  and 300 $\mu$m  wide, the laser propagates with an incidence angle of $19^{0}.$ The incoming boundary condition $\alpha=\alpha (y)$ is independent of time and mimics a laser beam whose width is equal to $ 40\mu m$ with five speckles ; each  speckle is modeled by a centered Gaussian function $h$  and is characterized by a random phase $\zeta_k$, that is to say $\alpha (y) =\Sigma_{k=1}^5 a_k  h(y-y_k) e^{i \zeta_k} ,$
 where the $\alpha_k$ are random  and the $a_k$ are close to each other. The plasma has an initial density equal to $N_m=0.15$ and the temperature is equal to $35. \ 10^6$ Kelvin. 
The mesh consists of 4 millions of cells and the time step is in the order of 0.1 picosecond (it is determined at each time step by the Courant-Friedrichs-Levy condition related to the sound speed of the plasma).
 The initial value of the laser intensity is zero, the plasma is progressively grabed by the ponderomotive force and on Figure \ref{beau}, we have plotted the laser intensity at different times. 
At the first snapshot (at time $2.6$ ps), the plasma is not grabed enough, so the value of $\mu$ is small; the autofocusing effect is very low but not negligible: instead of five different speckles at the incoming boundary  one notices only four speckels at the rear side (one of the four has a larger intensity) and a little spreading of the beam may be observed. At the second snapshot, the position of the four speckles has changed and the plasma is more grabed - since the energy density is larger in one speckle.
On the two last snapshots, we may check that the spreading of the beam at the rear side of the simulation box becomes larger when the time increases. Moreover the configuration is not stationary, this situation is characteristic of the so-called {\it filamentation instability.}

\section*{ Conclusion}
A mathematical analysis  has lead to an analytical form of the solution of the tilted paraxial equation in the simple case where the refraction index and the absorption coefficients are constant. Afterwards, we proposed a numerical method for solving the initial problem which uses the previous analytical form. 
The scheme has the property to yield a classical scheme when incidence angle becomes zero and the equation reduces to the classical paraxial one.
 The numerical method is illustrated by some results on toy problems. We have also given extensions of this model, which have enlarged the capability of our plateform {\it HERA} for laser propagation in a plasma (see \cite{fcs} and \cite{feu} for examples of simulations performed with {\it HERA}). This numerical method may be also extended in the case where the unit vector $\mathbf{K}$ depends slowly on the one-dimension spatial variable $\mathbf{x}.\mathbf{n}$, for instance if one has to deal with an equation of the following type
\begin{equation*}
 i\mathbf{K}.\nabla u +  i \frac1{2} (\nabla .\f{K}) u
+ \frac{1}{2k_0} \Delta ^{\f{k}} _{\bot }u- \mu u+i\nu u =0,  \quad   {\rm on} \;  {\cal D}.
\end{equation*}
The paraxial equation in a tilted frame may be also considered in a first region where the plasma density is slowly varying with respect to the spatial variable and coupled with another model in a neighbor region  where the plasma density is strongly varying: in that region the laser is no more characterized by the time-space envelope of the fast oscillating electric field but by the wave equation (\ref{wave})  (see \cite{desr}, for results obtained in {\it HERA} with this model). For simulating such a physical tilted beam, a classical paraxial model without accounting for the incidence angle would lead to search a the solution which would be highly oscillating with respect to the space variable and therefore to increase dramatically the mesh size to get accurate results.

\end{document}